\documentclass[a4paper,10pt]{amsart}

\usepackage[centertags]{amsmath}
\usepackage{amsfonts}
\usepackage{amssymb}
\usepackage{amsthm}
\usepackage{newlfont}
\usepackage{euscript}

\newcommand{\1}{1\!\!{\mathrm I}}
\renewcommand{\Re}{{\Bbb R}}
\newcommand{\eps}{\varepsilon}

\newcommand{\ax}{\Re^+}
\newcommand{\prt}{\partial}

\newcommand{\Ff}{{\EuScript F}}
\newcommand{\Ef}{{\EuScript E}}
\newcommand{\Bf}{{\EuScript B}}
\newcommand{\Tf}{{\EuScript T}}
\newcommand{\Cf}{{\EuScript C}}
\newcommand{\Yf}{{\EuScript Y}}

\newcommand{\Df}{{\EuScript D}}
\newcommand{\Nf}{{\EuScript N}}
\newcommand{\Gf}{{\EuScript G}}
\newcommand{\Qf}{{\EuScript Q}}
\newcommand{\Lf}{{\EuScript L}}
\newcommand{\Uf}{{\EuScript U}}
\newcommand{\Sf}{{\EuScript S}}
\newcommand{\Of}{{\EuScript O}}
\newcommand{\pf}{\UU_{fin}}

\newcommand{\ZZ}{{\Bbb Z}}
\newcommand{\NN}{{\Bbb N}}
\newcommand{\QQ}{{\Bbb Q}}
\newcommand{\TT}{{\Bbb T}}

\newcommand{\DD}{{\Bbb D}}
\newcommand{\UU}{{\Bbb U}}
\newcommand{\VV}{{\Bbb V}}
\renewcommand{\SS}{{\Bbb S}}
\newcommand{\ups}{\upsilon}
\newcommand{\Ups}{\Upsilon}

\newcommand{\demo}{\emph{Proof.} }
\newcommand{\Cd}{\mathsf{C}\,}
\newcommand{\Cs}{\mathsf{C}\,}
\newcommand{\Ls}{\mathsf{L}\,}
\newcommand{\Span}{\mathrm{span}\,}

\newcommand{\be}{\begin{equation}}
\newcommand{\ee}{\end{equation}}
\newcommand{\eqdef}{\mathop{=}\limits^{df}}
\newcommand{\tov}{\mathop{\to}\limits^{var}}
\renewcommand{\emptyset}{\varnothing}

\begin{document}

\numberwithin{equation}{section}
 \theoremstyle{plain}
\newtheorem{thm}{Theorem}[section]
\newtheorem{cor}{Corollary}[section]
\newtheorem{lem}{Lemma}[section]
\newtheorem{prop}{Proposition}[section]
\theoremstyle{definition}
\newtheorem{dfn}{Definition}[section]
\newtheorem{ex}{Example}[section]
\theoremstyle{remark}
\newtheorem{rem}{Remark}[section]

\title[Absolute continuity and convergence in variation] {Absolute continuity and
convergence in variation for distributions of a functionals of
Poisson point measure}
\author{Alexey M.Kulik}%
\address{Institute of Mathematics,
Ukrai\-ni\-an National Academy of Sciences, 3, Tereshchenkivska Str., Kyiv 01601, Ukraine}
 \abstract{General
 sufficient conditions are given for absolute
continuity and convergence in variation of distributions of a
functionals on a probability space, generated by a Poisson point
measure. The phase space of the Poisson point measure is supposed
to be of the form $\ax\times \UU$, and its intensity measure to be
equal $dt\Pi(du)$. We introduce the family of \emph{time
stretching transformations} of the configurations of the point
measure. The sufficient conditions for absolute continuity and
convergence in variation are given in the terms of the time
stretching transformations and the relative differential
operators. These conditions are applied to solutions of SDE's
driven by  Poisson point measures, including an SDE's with
non-constant jump rate. }
\endabstract
\email{kulik@imath.kiev.ua}%
 \maketitle

\section{Introduction} In this paper, we give a general
and transparent sufficient conditions for absolute continuity and
convergence in variation of a distributions of a functionals on
the probability space, generated by a Poisson point measure. The
phase space of the Poisson point measure is supposed to be of the
form $\ax\times \UU$, and its intensity measure to be equal
$dt\Pi(du)$, with $(\UU,\Uf)$ being Borel measurable space and
$\Pi$ being a $\sigma$-finite measure on $\Uf$. The Poisson point
measures of such a type arise naturally when the Levy processes or
their various modifications are considered; typically,
$\UU=\Re^m\backslash \{0\}$. For the Poisson point measures of
such a type, we introduce the family of \emph{time stretching
transformations}. The sufficient conditions for absolute
continuity and convergence in variation are given in the terms of
the time stretching transformations and the relative differential
operators. We illustrate the sufficient conditions obtained in
this paper by applying them to solutions of SDE's driven by
Poisson point measures, including SDE's with non-constant jump
rate.

Our approach strongly relies on an appropriate modification of
Yu.Davydov's \emph{stratification method}. This method is based on
disintegration of the probability space and finite-dimensional
change-of-variables formula. It is known that the stratification
method, unlike  the \emph{Malliavin calculus}, does not allow one
to prove the distribution density to be bounded, smooth, etc. The
main advantage of the stratification method is that it can be
applied under a very mild differentiability conditions on the
functionals under investigation, while, as we will see below, the
differential properties of the functionals of the Poisson point
measure typically are rather poor. In addition, this method
appears to be powerful enough to provide not only absolute
continuity for an individual distribution, but also convergence in
variation for a sequence of distributions. The latter finds a very
natural and useful applications in ergodic theory for  SDE's with
jump noise. See \cite{Me_ergodic}, where  the time stretching
transformations and associated stratifications are used to provide
the \emph{local Doeblin condition} for the solution to an SDE with
jumps, considered as a Markov process, and then to establish
ergodic and mixing rates for this process.

 This paper unifies and generalizes the
previous papers \cite{Me_1996},\cite{Me_TViMc},\cite{Me_conv_var}
by the same author. It contains partially the unpublished preprint
\cite{Me_preprint}.  The paper is also closely related to the
papers  \cite{Denis} and \cite{Nou_sim}. Statements 1 and 3 of
Theorem \ref{t933} below contain  Theorem A \cite{Nou_sim} as a
partial case. Statement 1 of Theorem \ref{t934} below  is a
generalization of Theorem 3.3.2 \cite{Denis}. However, Theorem
3.3.2 \cite{Denis} has a serious "gap" in its proof, that
seemingly can not be fixed up in the framework of \cite{Denis},
based on the Dirichlet form technique (see discussion in
subsection 4.3 below). The discussion of the relation between
Theorem \ref{t934} and the recent papers \cite{bally},
\cite{fourn_2008}, devoted to investigation of the laws of
solutions to SDE's with non-constant jump rate, is given in
subsection 4.2 below.

Let us give a brief overview of the other references related to
our investigation. The integration-by-parts structure for the pure
Poisson process was introduced independently in \cite{Carl_Pard}
and \cite{TsoiEl}. One can say that this structure, as well as its
extension used in our considerations, is provided by the \emph{
time-wise regularity} of the Poisson point measure. When the
Poisson point measure possesses some kind of a \emph{spatial
regularity} w.r.t. component $u\in\UU$,  other methods for
studying the local properties of the distributions of the
functionals are available, based both on the Malliavin-type
calculus and on the stratification technique. For exact
formulations, detailed discussion and further references at the
field, we refer to
\cite{bict_grav_Jac},\cite{bismut},\cite{Dav_Lif},\cite{Dav_Lif_Smor},\cite{kom_takeuchi},
\cite{leandre}. We also mention the method, introduced  by
J.Picard in \cite{picard} (see also
\cite{ishikawa},\cite{ishikawa_Kunita}), that,  in the case
$\UU=\Re^m\backslash\{0\}$, allows the L\'evy measure $\Pi$  of
the Poisson point measure to be singular, but requires some kind
of a \emph{frequency regularity} at the vicinity of $0$.

\section{Basic constructions and the main results}

\subsection{Basic constructions.}
 Everywhere below, $\UU$ is supposed to be a
locally compact metric space and $\Pi$ to be a $\sigma$-finite
measure on $\Bf(\UU)$, being finite on every bounded
$U\in\Bf(\UU)$. These suppositions,  if to compare with those made
in the Introduction, do not restrict generality, since one can
reduce Borel measurable space $(\UU,\Uf)$ with a $\sigma$-finite
measure $\Pi$ to $((0,1),\Bf(0,1))$ with a locally finite measure
$\Pi'$ by an appropriate Borel isomorphizm.

By $\nu$, we denote the Poisson point measure on $\ax\times\UU$
with its intensity measure equal $dt\Pi(du)$. By
$\Of\eqdef\Of(\ax\times\UU)$, we denote the space of
\emph{configurations} over $\ax\times\UU$, i.e. a family of
locally finite subsets of $\ax\times\UU$. The space $\Of$ is
equipped with the \emph{vague topology}, i.e. the weakest topology
such that every function $\Of\ni\varpi \mapsto \sum_{(t,u)\in
\varpi}f(t,u)$ with $f:\ax\times\UU\to \Re$ being a continuous
function with bounded support, is continuous. We denote $\Bf(\Of)$
the Borel $\sigma$-algebra on $\Of$ and write $P_\nu$ for the
distribution of the random element in $(\Of,\Bf(\Of))$ generated
by $\nu$. For more details, see e.g. \cite{Kerstan_Mattes_Mecke}.
  In the sequel, we suppose the basic probability space to have  the form
$(\Omega,\Ff,P)=(\Of,\Bf(\Of),P_\nu)$ and put
$\nu(\omega)=\omega$.

Denote $H=L_2(\Re^+), H_0=L_\infty(\Re^+)\cap L_2(\Re^+),
Jh(\cdot)=\int_0^\cdot h(s)\,ds,h\in H.$ For a fixed $h\in H_0$, we
define  the family $\{T_h^t,t\in\Re\}$ of transformations of the
axis $\Re^+$ by putting $T^t_hx, x\in\ax$ equal to the value at
the point $s=t$ of the solution of the Cauchy problem
\be\label{21} z'_{x,h}(s)=Jh(z_{x,h}(s)),\quad s\in \Re, \qquad
z_{x,h}(0)=x.\ee
 Since (\ref{21}) is the Cauchy problem for the
time-homogeneous ODE, one has that $T^{s+t}_h=T^s_h\circ T^t_h$,
and in particular $T_h^{-t}$ is the inverse transformation to
$T_h^t$. By multiplying $h$ by some $a>0$, we multiply, in fact,
the symbol $Jh(\cdot)$ of the equation by $a$. Now, making the
time change $\tilde s={s\over a}$, we see that
$T_h^a=T_{ah}^1,a>0$, which together with the previous
considerations gives that $T_h^t=T_{th}^1, h\in H_0, t\in\Re$.

Denote $T_h\equiv T_h^1$, we have just demonstrated that
$T_{sh}\circ T_{th}=T_{(s+t)h}$. This means that $\Tf_h\equiv
\{T_{th}, t\in \Re\}$ is a one-dimensional group of
transformations of the time axis $\ax$. It follows from the
construction that \be\label{211}{d\over dt}
|_{t=0}T_{th}x=Jh(x),\quad  x\in\ax.\ee

\begin{rem} We call $T_h$ the \emph{time stretching transformation}
because, for $h\in C(\ax)\cap H_0$, it can be constructed in a
more illustrative way: take the sequence of partitions $\{S^n\}$
of $\ax$ with $|S_n|\to 0, n\to +\infty$. For every $n$, we make
the following transformation of the axis: while preserving an
initial order of the segments, every segment of the partition
should be stretched by $e^{h(\theta)}$ times, where $\theta$ is
some inner point of the segment (if $h(\theta)<0$ then the segment
is in fact contracted). After passing to the limit (the formal
proof is omitted here in order to shorten the exposition) we
obtain the transformation $T_h.$ Thus one can say that $T_h$
performs the stretching of every infinitesimal segment $dx$ by
$e^{h(x)}$ times.
\end{rem}

 Denote $\UU_{fin}=
 \{\Gamma\in \Bf(\UU),\Gamma$ is bounded$\}$ and
define, for $h\in H_0, \Gamma\in \pf$, a transformation $T_h^\Gamma$
of the random measure $\nu$ by
$$
[T_h^\Gamma \nu]([0,t]\times \Delta)=\nu([0,T_{h}t]\times
(\Delta\cap\Gamma))+ \nu([0,t]\times (\Delta\backslash\Gamma))
,\quad t\in\Re^+,\Delta\in \pf.
$$
This transformation is  generated by the following transformation
of the space of the configuration: $(\tau,x)\in\omega$ with
$x\not\in \Gamma$ remains unchanged; for every point $(\tau,x)\in
\omega$ with $x\in \Gamma$, its ``moment of the jump" $\tau$ is
transformed to $T_{-h}\tau$; neither any point of the
configuration is eliminated nor any new point is added to the
configuration. In the sequel we denote, by the same symbol
$T_h^\Gamma$, the bijective transformation of the space of
configurations described above.

The image $T_h^\Gamma \nu$ is again a random Poisson point
measure, and its intensity measure can be expressed through $\Pi$
and $r_h(x)\eqdef {d\over dx}(T_hx)$ explicitly. An easy
calculation gives that
 \be\label{212}r_h(x)=
 \int_0^1 h(T_{sh}x)\,ds,\quad x\in\Re^+.\ee
Thus the following statement is a corollary of the classical
absolute continuity result for L\'evy processes, see
\cite{Skor_nez_prir}, Chapter 9. We put
 $$
  p_h^\Gamma=\exp\left\{\int_{\Re^+}
 r_h(t)\nu(dt,\Gamma)-\lim_{t\to+\infty}[T_ht-t]\Pi(\Gamma)\right\}.
 $$

\begin{lem}\label{l21} The transformation $T_h^\Gamma$ is admissible for
the distribution of $\nu$ with the density $p_h^\Gamma$, i.e., for
every $\{t_1,\dots,t_n\}\subset\ax,
\{\Delta_1,\dots,\Delta_n\}\subset \pf$ and  Borel function
$\phi:\Re^n\to \Re$,
$$
E\phi( [T_h^\Gamma \nu]([0,t_1]\times \Delta_1),\dots,[T_h^\Gamma
\nu]([0,t_n]\times \Delta_n))=Ep_h^\Gamma\phi(\nu ([0,t_1]\times
\Delta_1),\dots,\nu([0,t_n]\times \Delta_n)).
$$
\end{lem}

The lemma implies that every transformation $T_h^\Gamma$ generates
the corresponding transformation of the random variables. In the
sequel, we denote the latter transformation by the same symbol
$T_h^\Gamma$.

\begin{dfn}\label{d21} Let $h\in H_0, \Gamma\in \pf$ be fixed.

 1. The functional $f\in L_0(\Omega,\Ff,P)$ is said to be \emph{almost
surely} (\emph{a.s.}) \emph{differentiable} in the direction
$(h,\Gamma)$ and to have  almost sure (a.s.) derivative
$D_h^\Gamma f$, if
 \be\label{213} {T_{\eps h}^\Gamma f-f\over
\eps}\to D_h^\Gamma f,\quad \eps\to 0 \ee almost surely.

2. Let $p\in [1,+\infty)$. The functional $f\in L_p(\Omega,\Ff,P)$
is said to be $L_p$-\emph{differentiable} in the direction
$(h,\Gamma)$ and to have $L_p$ derivative $D_h^\Gamma f$, if
convergence (\ref{213}) holds in $L_p$ sense.
\end{dfn}

Let us give an example demonstrating one specific property of the
family $\{T_h, h\in H_0\}$.

\begin{ex}\label{e25} Let $f=\tau_n^\Gamma\eqdef\inf\{t:\nu(t,\Gamma)=n\}$, and let
 $h,g\in C(\ax)\cap
H_0$ be such that $h(t)\int_0^tg(s)\,ds\not=g(t)\int_0^t g(s)\,
ds, t>0$. Then
$$
D^\Gamma_h D^\Gamma_g f=h(\tau^\Gamma_n)\int_0^{\tau_n^\Gamma}
g(s)\,ds\not= g(\tau^\Gamma_n)\int_0^{\tau_n^\Gamma}
h(s)\,ds=D^\Gamma_g D^\Gamma_h f
$$
almost surely (this follows from the relation (\ref{221}) given
below). In particular, the family of transformations
$\{T^\Gamma_h, h\in H_0\}$ is not commutative and therefore cannot
be considered as an infinite-dimensional additive group of
transformations.
\end{ex}

The non-commutative structure of the family $\{T_h, h\in H_0\}$
does not allow one to apply the stratification method for study of
the absolute continuity of the laws of differentiable functionals
straightforwardly. In order to overcome this difficulty, we
introduce an additional construction based on the notion of a
\emph{differential grid}.

\begin{dfn}\label{d26}
 A family $\Gf=\{[a_i,b_i)\subset \ax, h_i\in
H_0, \Gamma_{i}\in\pf, i\leq m\}$ is called  \emph{a differential
grid} (or simply \emph{a grid}) if

(i) for every $i\not= j$,  $\Bigl([a_i,b_i)\times \Gamma_i\Bigr)
\cap \Bigl( [a_j,b_j)\times \Gamma_j\Bigr) =\emptyset$;

(ii) for every $i\in\NN$, $Jh_i>0$ inside $(a_i,b_i)$ and $Jh_i=0$
outside $(a_i,b_i)$.

\noindent The number $m\in \NN$ is called a dimension of the grid
$\Gf$.
\end{dfn}

Denote $T_t^{i}=T_{th_{i}}^{\Gamma_{i}}$.  It follows from the
construction of the transformations $T_h^\Gamma$ that, for a given
$i\in \NN, t\in \Re$,
$$
T_t^{i}\tau_n^{\Gamma_{i}}=T_{t h_i}\tau_n^{\Gamma_{i}}\quad
\begin{cases}
=\tau_n^{\Gamma_{i}},& \tau_n^{\Gamma_{i}}\not\in [a_i,b_i)\\
\in [a_i,b_i),& \tau_n^{\Gamma_{i}}\in [a_i,b_i)
\end{cases}\quad \hbox{for every } n.
$$
In other words: a grid $\Gf$ generates a partition of some part of
the phase space $\ax\times \UU$ of the random measure $\nu$ into
the non-intersecting  cells $\{\Gf_{i}=[a_i,b_i)\times
\Gamma_{i}\}$. The transformation $T_t^{i}$ does not change points
of configuration outside the cell $\Gf_{i}$ and keeps the points
from this cell in it. In addition, for every $i\leq m, t,\tilde
t\in \Re $, the transformations $T_t^{i}$,$T_{\tilde t}^{ i}$
commute because so do the time axis transformations $T_{t
h_i}$,$T_{\tilde t h_i}$. Therefore, for every $i,\tilde i \leq m,
t,\tilde t\in \Re $, the transformations $T_t^{i}$,$T_{\tilde
t}^{\tilde i}$ commute. This implies the following proposition.

\begin{prop}\label{p270} For a given grid $\Gf$ and $z=(z_{i})_{i=1,\dots,m}\in \Re^m$,
define the transformation
$$
T^{\Gf}_{z}=T^{1}_{z_{1}}\circ T^{2}_{z_{2}}\circ\dots \circ
T^{m}_{z_{m}}.
$$
Then $\Tf^{\Gf}=\{T^\Gf_z, z\in\Re^m\}$ is the group of admissible
transformations of $\Omega$ which is additive in the sense that
$T^\Gf_{z^1+z^2}=T^\Gf_{z^1}\circ T^\Gf_{z^2}, z^{1,2}\in\Re^m.$
\end{prop}

It can be said that, by fixing the grid $\Gf$, we choose from the
whole variety of admissible transformations $\{T_h^\Gamma, h\in
H_0,\Gamma\in \pf\}$ the additive subfamily that is more
convenient to deal with.

\begin{dfn}\label{d28} 1. The functional $f\in L_0(\Omega,\Ff,P)$ is
a.s. stochastically differentiable w.r.t. differential grid $\Gf$
if $f$ is a.s. differentiable in every direction $(h_i, \Gamma_i),
i=1,\dots, m$.  The random vector $D^\Gf
f=(D^{\Gamma_1}_{h_1}f,\dots, D^{\Gamma_m}_{h_m}f)$ is called the
a.s. stochastic derivative of $f$.

2. Let $p\in[1,+\infty)$.  The functional $f\in L_p(\Omega,\Ff,P)$
is stochastically differentiable in $L_p$ sense w.r.t.
differential grid $\Gf$ if $f$ is $L_p$-differentiable in every
direction $(h_i, \Gamma_i), i=1,\dots, m$.  The random vector
$D^\Gf f=(D^{\Gamma_1}_{h_1}f,\dots, D^{\Gamma_m}_{h_m}f)$ is
called the $L_p$ stochastic derivative of $f$.
\end{dfn}

We denote $D_i^\Gf f=D_{h_i}^{\Gamma_i} f,i=1,\dots,m$.

\subsection{Sufficient conditions for absolute continuity and convergence in variation.}

The proofs for the following  theorems are given in Section 3
below.

\begin{thm}\label{t212}
Consider an $\Re^m$-valued  random vector $f=(f_1,\dots,f_m)$ and
a grid $\Gf$ of dimension $m$. Let every component of the vector
$f$ to be differentiable w.r.t. $\Gf$ either in a.s. or in $L_p$
sense for some $p\geq 1$.

Denote $\Sigma^{f,\Gf}=(D^\Gf_i f_j)_{i,j=1}^m$ and put $
\Nf(f,\Gf)=\{\omega:\hbox{ the matrix
}\Sigma^{f,\Gf}(\omega)\hbox{ is non-degenerate}\}$. Then
$$
P\Bigl|_{\Nf(f,\Gf)}\circ f^{-1}\ll \lambda^m.
$$
\end{thm}

\begin{thm}\label{t217}
Consider a sequence of $\Re^m$-valued random vectors $\{f^n, n\geq
1\}$ such that, for a given grid $\Gf$ of dimension $m$, every
component $f_j^n, j=1,\dots, m$ of the vector $f^n$ is $L_m$
differentiable w.r.t. $\Gf$. Suppose that
$$
f^n_j\to f_j,\quad D_i^\Gf f^n_j\to D_i^\Gf f_j\hbox{ in
}L_{m},\quad n\to +\infty,\quad i,j=1,\dots,m.
$$

 Then, for every $A\subset{\Nf(f,\Gf)}$,
$$
P\Bigl|_{A}\circ f^{-1}_n\to P\Bigl|_{A}\circ f^{-1},\quad n\to
+\infty
$$
in variation.
\end{thm}

We remark that the type of differentiability of the components of
$f$ is unimportant in the condition for absolute continuity, given
in Theorem \ref{t212}. On the contrary, this type is crucial in
the condition for convergence in variation. For instance, the
immediate analogue of Theorem \ref{t217}, with the $L_m$
derivatives replaced by the a.s. ones, fails to be true. One can
construct the counterexample to such a statement using Example 1.2
\cite{Me_conv_var}. In order to formulate the correct version of
Theorem \ref{t217} in the terms of  a.s. derivatives, we need an
auxiliary notion.

\begin{dfn}\label{d31} The sequence of the measurable functions
$\{f_n:\Omega\to \Re, n\geq 1\}$ is said to have \emph{a uniformly
dominated increments} w.r.t. the grid $\Gf$ on the set
$\Omega'\in\Ff$, if there exist a random variable $\varrho$ and  a
family of jointly measurable functions $\{g_{i}:\Omega\times
\Re\to \Re\}$ such that

(i) for every $i$ and almost every $\omega$, the function
$g_{i}(\omega,\cdot)$ is an increasing one;

(ii) $\varrho>0$ almost surely and,  for every $n\geq
1,\omega\in\Omega $, \be\label{938} |T_{t}^{\Gf}
f_n(\omega)-T_{s}^{\Gf} f_n(\omega)|\leq
\sum_{i=1}^m\Big[g_{i}(\omega, t_i\vee s_i)-g_{i}(\omega,
t_i\wedge s_i)\Big], \quad \|t\|,\|s\|<\varrho(\omega),
\,T_{t}^{\Gf}\omega\in\Omega', \,T_{t}^{\Gf}\omega\in\Omega'. \ee
A sequence of $\Re^m$-valued random vectors $\{f^n, n\geq 1\}$ is
said to have a uniformly dominated increments w.r.t. the grid
$\Gf$ on the set $\Omega'\in\Ff$ if every sequence $\{f^n_j, n\geq
1\}, j=1,\dots,m$ has a uniformly dominated increments w.r.t. the
grid $\Gf$ on this set.
\end{dfn}

\begin{thm}\label{t218}
Consider a sequence of $\Re^m$-valued random vectors $\{f^n, n\geq
1\}$ such that, for a given grid $\Gf$ of dimension $m$, every
component $f_j^n, j=1,\dots, m$ of the vector $f^n$ is a.s.
differentiable w.r.t. $\Gf$. Suppose that
$$
f^n_j\to f_j,\quad D_i^\Gf f^n_j\to D_i^\Gf f_j\hbox{ in
probability},\quad n\to +\infty,\quad i,j=1,\dots,m.
$$
Suppose additionally that $\{f^n\}$ has a uniformly dominated
increments w.r.t. $\Gf$ on the set $\Omega'$.

Then, for every $A\subset{\Nf(f,\Gf)}\cap \Omega'$,
$$
P\Bigl|_{A}\circ f^{-1}_n\to P\Bigl|_{A}\circ f^{-1},\quad n\to
+\infty
$$
in variation.
\end{thm}

\section{The stratification method and proofs of Theorems \ref{t212} -- \ref{t218}}

In this section, we prove the general statements, formulated in
section 2.2. Our main tool is a certain version of Yu.Davydov's
stratification method; for the basic constructions of this method,
references and further discussion, we refer the reader to the
Chapter 2 of the monograph \cite{Dav_Lif_Smor}. Some steps in our
considerations have an analogues in the available literature. For
instance, the trick with using Theorem 3.1.16 \cite{Federer} in
order to replace a function differentiable in some weak sense by a
$C^1$ one, was used in \cite{Boga_Smol} in the context of
stratifications generated by linear shifts and in
\cite{Bouleau_hirsch}, Chapter II.5 in the context of Dirichlet
forms on vector spaces.

\subsection{Stratifications, generated by differential grids.}

Let $\Gf$ be a differential grid of dimension $m$. For every
$\omega\in \Omega$, consider the set
$\upsilon=\upsilon(\omega)=\{T_z^\Gf\omega, z\in \Re^m\}$. This
set is called \emph{the orbit} of the group $\Tf^\Gf$,
corresponding to $\omega$. The set of all such an orbits is
denoted $\Upsilon$, and $\Omega$ is represented as the disjunctive
union \be\label{931}\Omega=\bigsqcup_{\ups\in\Ups}\ups.\ee The
decomposition (\ref{931}) is called the \emph{stratification} of
$\Omega$ to the orbits of the group $\Tf^\Gf$.

Every orbit $\ups$ has a simple structure. Denote
$\Df^\Gamma\eqdef\{\tau\in\Df: p(\tau)\in\Gamma\}$,
$$
I(\omega)=\Big\{i=1,\dots, m:
(a_i,b_i)\cap\Df^{\Gamma_i}(\omega)\not=\varnothing\Big\},\quad
\omega\in \Omega.
$$
By condition (ii) of Definition \ref{d26}, for every
$i=1,\dots,m$, the mapping
$$
\Re\ni t\mapsto T_{th_i}x\in \ax
$$
is the identical one if $x\not\in (a_i,b_i)$ and is strictly
monotonous if $x\in (a_i,b_i)$. This implies the following
equivalence: for every $z^1,z^2\in \Re^m, \omega\in \Omega$,
$$
T_{z^1}^\Gf\omega=T_{z_2}^\Gf\omega\Longleftrightarrow
z^1_i=z^2_i,\quad i\in I(\omega).
$$
Therefore, the orbit $\ups(\omega)$ is the bijective image of
$\Re^{\#I(\omega)}$ (here and below, $\#$ is used for the number
of elements of the set).

Denote $\Omega^\Gf=\Big\{\omega: I(\omega)=\{1,\dots, m\}\Big\}$,
one can see that $\Omega^\Gf$ is measurable. For our further
purposes, it would be enough to restrict the initial probability
$P$ to $\Omega^\Gf$ and to describe the stratification of
$\Omega^\Gf$, only. Such a restriction simplifies the exposition,
since, for every point $\omega \in\Omega^\Gf$, the corresponding
orbit is a bijective image of $\Re^{m}$.

\begin{lem}\label{l31} There exists a complete separable metric
space $\Yf$ and a bijection $\vartheta:\Omega\to \Yf\times \Re^m$
such that $\vartheta$ is $\Ff$ -- $\Bf(\Yf)\otimes \Bf(\Re^m)$
measurable, $\vartheta^{-1}$ is $\Bf(\Yf)\otimes \Bf(\Re^m)$ --
$\Ff$ measurable and \be\label{932}
\vartheta\Big(T_{z}^\Gf\omega\Big)=\Big(\pi_1(\vartheta(\omega)),\pi_2(\vartheta(\omega))+z
\Big),\quad \omega\in \Omega^\Gf, \, z\in\Re^m,\ee where
$\pi_1,\pi_2$ denote the projections on the first and the second
coordinates in $\Yf\times \Re^m$ respectively.
\end{lem}

\demo First of all, we mention that $\Omega=\Of$ can be considered
as a Polish space via the following construction.  For two
configurations $\omega',\omega''\in\Of$, we put
$$
d_\Of(\omega',\omega'')=\sum_{m=1}^\infty 2^{-m}[1\wedge
d_H(\omega'\cap K_m,\omega''\cap K_m)],
$$
where $d_H$ is the Haussdorff metrics on the set of closed subsets
of $\ax\times \UU$, and $\{K_m\}$ is a sequence of compacts such
that $\bigcup_m K_m=\ax\times \UU$. Then $(\Of,d_\Of)$ is a Polish
space, and one can deduce from \cite{materon}, Propositions 1.4.1
and 1.4.4, that the Borel structure  on $\Of$ generated by $d_\Of$
coincides with the one generated by the vague topology.

 For $\omega\in\Omega^\Gf$, define
$\tau_i(\omega)=\min\Big[(a_i,b_i)\cap\Df^{\Gamma_i}(\omega)\Big],
i=1,\dots,m.$ We have already mentioned that,  for every
$i=1,\dots,m$ and $x\in(a_i,b_i)$, the transformation $\Re\ni
t\mapsto T_{th_i}x \in (a_i,b_i)$ is strictly monotonous. In
addition, it is bijective and continuous together with its
inverse. Therefore, for every $i=1,\dots,m$ there exists unique
$z_i(\omega)\in\Re$ such that
$T_{-z_i(\omega)h_i}^{\Gamma_i}\tau_i(\omega)={1\over
2}(a_i+b_i)$. Denote
$z(\omega)=(z_1(\omega),\dots,z_m(\omega))\in\Re^m$. Denote by
$\Yf$ the  family of all configurations satisfying the following
additional condition: for every set $(a_i,b_i)\times \Gamma_i,
i=1,\dots,m$, the configuration is not empty in this set, and the
smallest time coordinate of the point in this set is equal to
${1\over 2}(a_i+b_i)$. This family is a complete separable metric
space w.r.t. the local Haussdorf metrics described above.

Now put
$$\vartheta(\omega)=\Big(
T_{-z(\omega)}^\Gf\omega,z(\omega)\Big)\in \Yf\times \Re^m, \quad
\omega\in\Omega^\Gf.
$$
The map $\vartheta$ is a bijection between $\Omega^\Gf$ and
$\Yf\times \Re^m$. One can easily see that both $\vartheta$ and
$\vartheta^{-1}$ are measurable (moreover, continuous). At last,
if $\tilde \omega=T^\Gf_{\tilde z}\omega$ then, by the group
property of the family $\Tf^\Gf$,
$$
z(\tilde \omega)=z(\omega)+\tilde z.
$$
This proves (\ref{932}). The lemma is proved.

In a sequel, we denote the points of $\Yf$  by $\ups$ in order to
emphasize that $\Yf$, in fact,  is the set of the orbits. We also
omit $\vartheta$  in the notation and identify $\omega\in
\Omega^\Gf$ with its image $(\ups,z)\in \Yf\times \Re^m$.

For $A\subset \Omega^\Gf,\ups\in\Yf$ denote $A_\ups=\{z\in \Re^m:
(\ups,z)\in A\}$. Similarly, for $\ups\in\Yf$ and the function
$f:\Omega^\Gf\to \Re$, define the function $f_\ups: \Re^m\ni
l\mapsto f((\ups,z))\in\Re^m$. It follows from (\ref{932}) that
\be\label{934} [T_z^\Gf f]_\ups(\cdot)=f_\ups(\cdot+z),\quad
\ups\in\Yf, \, z\in \Re^m.\ee

Denote, by $P^\Gf$, both $P|_{\Omega^\Gf}$ and its image under
$\vartheta$. Denote, by $P^\Yf$, the projection of $P^\Gf$ on the
first coordinate in $\Yf\times \Re^m$ (i.e., the image of $P^\Gf$
under the projection $\pi_1$). The following statement is a
version of the well known theorem on existence of the family of
conditional distribution (e.g. \cite{Partas}, Chapter 5).

\begin{prop} There exists a family $\{P_\ups,\ups\in\Yf\}$ of finite
measures on $\Bf(\Re^m)$ such that

1. for every $B\in \Bf(\Re^m)$, the function $\ups\mapsto
P_\ups(B)$ is Borel measurable;

2. for every $A\in \Ff, A\subset \Omega^\Gf$,
$$
P(A)=\int_\Yf P_{\ups}(A_\ups)P^\Yf(d\ups).
$$
\end{prop}

\begin{lem}\label{l32} 1. For $P^\Yf$ -- almost all $\ups\in\Yf$,
the measure $P_\ups$ possesses a continuous strictly positive
density w.r.t. $\lambda^m$.

2. If $f$ is a.s. stochastically differentiable w.r.t. the grid
$\Gf$, then, for $P^\Yf$ -- almost all $\ups\in\Yf$, the function
$f_\ups$  $\lambda^m$-almost everywhere possesses partial
derivatives ${\prt\over \prt z_1}f_\ups,\dots,{\prt\over \prt
z_m}f_\ups$. In addition, \be\label{933}
 {\prt\over \prt
z_i}f_\ups=[D_{h_i}^{\Gamma_i} f]_\ups, \quad i=1,\dots,m\ee
almost surely.

3. If $f$ is stochastically differentiable in $L_p$ sense w.r.t.
the grid $\Gf$, then, for $P^\Yf$ -- almost all $\ups\in\Yf$, the
function $f_\ups$ belongs to the local Sobolev space
$W_{p,loc}^1(\Re^m,\lambda^m)$. In addition, the relation
(\ref{933}) holds true almost surely (in this case, ${\prt \over
\prt z_i}$ is the Sobolev partial derivative and
$D_{h_i}^{\Gamma_i}$ is the $L_p$ stochastic derivative).
\end{lem}

\demo Let $C^1_0(\Re^m)$ denote the set of continuously
differentiable functions $\Re^m\to \Re$  with a compact supports.
Denote, by $\Cf_0$, the set of measurable functions $g$ on
$\Omega$ such that $g_\ups\in C^1_0(\Re^m), \ups\in\Yf$ and
$$
\sup_{\ups\in\Yf}\sup_{z\in\Re^m} \Big[|g_\ups(z)|+\|\nabla
g_\ups(z)\|_{\Re^m}\Big]<+\infty.
$$
By (\ref{934}) and the dominated convergence theorem, every $g\in
\Cf_0$ is $L_p$ differentiable w.r.t. the grid $\Gf$ for every
$p\geq 1$, and (\ref{933}) holds true at every point.

Denote $\rho_i=-\int_{0}^\infty h_i(t)\tilde \nu(t,\Gamma_i),\,
i=1,\dots, m$. One can verify that
$$
{1-p_{\eps h_i}^{\Gamma_i}\over \eps}\to \rho_i,\quad \eps\to 0,
\quad i=1,\dots,m
$$
in $L_p$ sense for every $p\geq 1$. Therefore, for every $g\in
\Cf_0$, \be\label{935} E D_{h_i}^{\Gamma_i}g=\lim_{\eps\to
0}E{T_{\eps h_i}^{\Gamma_i}g-g\over \eps}=\lim_{\eps\to
0}Eg{p_{\eps h_i}^{\Gamma_i}-1\over \eps}=-Eg\rho_i,\quad
i=1,\dots, m. \ee

Consider a countable dense subset $\Phi$ of $C_0^1(\Re^m)$. For
every $\phi\in\Phi, C\in\Bf(\Yf)$, consider the function $\phi_C:
(\ups, z)\mapsto \1_C(\ups)\phi(z)$. This function belongs to
$\Cf_0$ and integration-by-parts formula (\ref{935}) for this
function has the form
$$
\int_C\int_{\Re^m}{\prt\over \prt
z_i}\phi(z)P_\ups(dz)P^\Yf(d\ups)=-\int_C\int_{\Re^m}[\rho_i]_\ups(z)
\phi(z)P_\ups(dz)P^\Yf(d\ups),\quad i=1,\dots,m.
$$
Since $C\in\Bf(\Yf)$ is arbitrary, we conclude that, for a given
$\phi\in \Phi$, \be\label{936} \int_{\Re^m}{\prt\over \prt
z_i}\phi(z)P_\ups(dz)P^\Yf(d\ups)=-\int_{\Re^m}[\rho_i]_\ups(z)
\phi(z)P_\ups(dz)P^\Yf(d\ups),\quad i=1,\dots,m\ee for
$P^\Yf$-almost all $\ups$. Denote, by $\Yf_\phi$, the set of
$\ups\in \Yf$ such that (\ref{936}) holds.

Every $\rho_i$ is an integral of a bounded function over a
compensated Poisson point measure of the finite intensity.
Therefore, $E\exp|\rho_i|<+\infty,i=1,\dots, m$. Then there exists
a set $\Yf_*$ with $P^\Yf(\Yf\backslash \Yf_*)=0$ such that
$$
\int_{\Re^m}\exp\Big|[\rho_i]_\ups(z)\Big|\, P_\ups(dz)<+\infty,
\quad i=1,\dots, m, \quad \ups\in \Yf_*.
$$
Thus, for $\ups\in \Yf^*\eqdef\Yf_*\cap \bigcap_{\phi}\Yf_\phi$,
the relation (\ref{936}) holds true for every $\phi\in
C_0^1(\Re^m)$ and every function $[\rho_i]_\ups$ possesses an
exponential moment. In other words: for every $\ups\in \Yf^*$, the
measure $P_\ups$ is \emph{differentiable} w.r.t. the basic
directions in $\Re^m$ and its \emph{logarithmic derivative}
possesses an exponential moment. Then Proposition 4.3.1
\cite{Boga} provides that $P_\ups$ possesses a continuous strictly
positive density. This completes the proof of the statement 1.
This statement provides that, for $P^\Yf$-almost all $\ups\in\Yf$,
 $P_\ups$-a.s. convergence is equivalent to
$\lambda^m$-a.s. convergence and  $L_p(\Re^m, P_\ups)$-convergence
implies  $L_{p,loc}(\Re^m, \lambda^m)$-convergence. Now the
statements 2,3 follow from (\ref{934}) and Fubini theorem. The
lemma is proved.

\subsection{Proof of Theorem \ref{t212}} For every $\omega\not\in\Omega^\Gf,$ there exists
$i\in\{1,\dots,m\}$ such that $T_{th_i}^{\Gamma_i}\omega=\omega,
t\in\Re$. This implies that, for such an $\omega$, at least one
column in the matrix $\Sigma^{f,\Gf}$ contains zeroes only. Thus,
$\Nf(f,\Gf)\subset \Omega^\Gf$ and \be\label{937} P(\{f\in B\}\cap
\Nf(f,\Gf))=\int_{\Yf}P_{\ups}(f_\ups\in B,\nabla f_\ups\hbox{ is
non-degenerate})P^{\Yf}(d\ups),\quad B\in \Bf(\Re^m), \ee here
$\nabla f_u$ denotes the matrix that contains the partial
derivatives of $f_u$, either a.s. or Sobolev ones. Here we have
used that, by Lemma \ref{l32},   $\nabla
f_\ups=[\Sigma^{f,\Gf}]_\ups$ almost surely on $\Re^m$ for
$P^\Yf$-almost all $\ups$.

 By  (\ref{937}) and Fubini theorem, it is enough to
prove that, for almost all $\ups$, the image under the mapping
$f_\ups$ of the Lebesgue measure restricted to
$\Nf(f_\ups)\eqdef\{x: \nabla f_\ups(x)$, is non-degenerate$\}$ is
absolutely continuous. The crucial step in the proof of the latter
fact is provided by the the following statement.
\begin{lem}\label{l933} Let $F:\Re^m\to\Re^n$ to have the approximative partial derivative
w.r.t. every basic direction  at $\lambda^m$-almost all points of
$x\in \Re^m$ and $G$ be the corresponding approximative gradient.
Then for every $\eps>0$ there exists  $F_\eps\in C^1(\Re^m,\Re^n)$
such that
$$ \lambda^m(\{x:
F(x)\not=F_\eps(x)\}\cup\{x: G(x)\not=\nabla F_\eps(x)\})<\eps.
$$
\end{lem}

Lemma \ref{l933} is a corollary of the following two statements,
given in \cite{Federer}.

\begin{prop}\label{p31} I. (\cite{Federer}, Theorem
3.1.4). Let the function $F:\Re^m\to \Re^n$  to possess the
approximative partial derivative w.r.t. every basic direction  at
all the points of a set $A\subset \Re^m$. Then, for
$\lambda^m$-almost all points $a\in A$, the function  $F$
possesses the approximative differential.

II. (\cite{Federer}, Theorem 3.1.16). Let  $A\subset \Re^m, f:A\to
\Re^n$ and
\be\label{51_ap_lim_sup}\mathop{ap\lim\sup}\limits_{x\to
a}{\|F(x)-F(a)\|_{\Re^n}\over \|x-a\|_{\Re^m}}<+\infty \ee for
$\lambda^m$-almost all $a\in A$, then, for every $\eps>0$, there
exists $F_\eps\in C^1(\Re^m, \Re^n)$ such that
$$ \lambda^m(\{x:F(x)\not=F_\eps(x)\}<\eps.
$$
\end{prop}

We do not discuss here the notions of the approximative upper
limit ($ap\lim\sup$), approximative partial derivative and
approximative differential, referring to \cite{Federer}, Chapter
3. We just mention that, if $F$ either belong to
$W_{p,loc}^1(\Re^m)$ or possesses partial derivatives w.r.t. basic
directions at $\lambda^m$-almost all points, then $F$ possesses
approximative partial derivatives w.r.t. basic directions at
$\lambda^m$-almost all points. Moreover, if the function  $F$
possesses the approximative differential at the point $a$, then
(\ref{51_ap_lim_sup}) holds true at this point.

Thus, for every $\eps>0$ and almost every $\ups$, there exists
$f^\eps_\ups\in C^1(\Re^m,\Re^m)$ such that the Lebesgue measure
of the  set
$$
C(\ups,\eps)\eqdef\{z\in \Re^m:f_\ups(z)\not=f_\ups^\eps(z)\}
$$
is less that $\eps$. By the Lebesgue theorem, almost every point
of $\Re^m\backslash C(\ups,\eps)$ is a Lebesgue point (i.e. a
density point), and therefore $\nabla f_\ups=\nabla f_{\ups}^\eps$
a.e. on $\Re^m\backslash C(\ups,\eps)$.

Now, the image measure of $\lambda^m|_{\Nf(f_\ups)}$ under
$f_\ups$ can be represented as the sum of the two measures
$$
\lambda_m|_{[\Re^m\backslash
C(\ups,\eps)]\cap\Nf(f_\ups^\eps)}\circ [f_\ups^\eps]^{-1}
\quad\hbox{ and
}\quad\lambda_m|_{C(\ups,\eps)\cap\Nf(f_\ups^\eps)}\circ
[f_\ups]^{-1}.
$$
The first one is absolutely continuous by the standard
change-of-variables formula for $C^1$-transformations. The second
one has its total mass being less than $\eps$. Since $\eps>0$ is
arbitrary, this proves the required absolute continuity. The
theorem is proved.

\subsection{Proofs of Theorems \ref{t217},\ref{t218}}
Theorem 2.1 \cite{Pilipenko_convergence_by_variation} provides the
criterium for convergence in variation of induced measures on a
finite-dimensional space. In our considerations, we use two
following sufficient conditions,  based on this criterium.

\begin{prop}\label{p27} I. (\cite{Pilipenko_convergence_by_variation}, Corollary
2.7). Let $F_,F_n\in W_{p,loc}^1(\Re^m, \Re^m)$ with $p\geq m$,
and  $F_n\to F,n\to\infty$ w.r.t. Sobolev norm
$\|\cdot\|_{W_p^1(\Re^m,\Re^m)}$ on every ball. Then \be\label{27}
\lambda^m|_A\circ F_n^{-1}\mathop{\longrightarrow}\limits^{var}
\lambda^m|_A\circ F^{-1}, n\to +\infty \hbox{ for every measurable
}A\subset \{\det \nabla F\not=0\}. \ee

II. (\cite{Me_conv_var}, Theorem 3.1). Let $F_n,F:\Re^m\to \Re^m$
possess approximative partial derivatives at $\lambda^m$-almost
every point and $F_n\to F,\nabla F_n\to \nabla F$ in a sense of
convergence in measure $\lambda^m$. Let, in addition, the sequence
$\{F_n\}$ be \textbf{uniformly approximatively Lipschitz}. This,
by definition, means that,  for every $\delta>0, R<+\infty,$ there
exist a compact set $K_{\delta,R}$ and a constant
$L_{\delta,R}<+\infty$ such that
$\lambda^m(B_{\Re^m}(0,R)\backslash K_\delta)<\delta$
 and every function $F_n|_{K_\delta}$ is a
Lipschitz function with the Lipschitz constant $L_{\delta,R}$.
Then (\ref{27}) holds true.
\end{prop}

Under conditions of Theorem \ref{t217}, the statements 1 and 3 of
Lemma \ref{l32} provide that, for $P^\Yf$-almost all $\ups$,
$f^n_\ups\to f_\ups$ w.r.t. Sobolev norm
$\|\cdot\|_{W_p^1(\Re^m,\Re^m)}$ on every ball. Then the statement
I of  Proposition \ref{p27} and the statement 1  of Lemma
\ref{l32} provide that, for every $B\in\Bf(\Re^m)$,
$$ P_{\ups}|_{\Nf(f_\ups)\cap B}\circ
[f^n]_\ups^{-1}\tov P_\ups|_{\Nf(f_\ups)\cap  B}\circ
[f^n]_\ups^{-1}\quad\hbox{for $P^\Yf$-almost all } \ups\in \Yf.$$
By applying the decomposition formula (\ref{937}) and Fubini
theorem, we complete the proof of Theorem  \ref{t217}.

Via the same arguments, the statement II of Proposition \ref{p27}
would provide the proof of Theorem \ref{t218}, but we have to
verify additionally that, for $P^\Yf$-almost all $\ups$, the
sequence $\{[f^n\1_{\Omega'}]_\ups\}$ is uniformly approximatively
Lipschitz.

Recall that $\Nf(f,\Gf)\subset \Omega^\Gf$ and thus we can exclude
$\omega\not\in \Omega^\Gf$ from the consideration. By the analogy
with Definition \ref{d31}, we say that the sequence of measurable
functions $\{F_n:\Re^m\to \Re, n\geq 1\}$ has \emph{a uniformly
dominated increments} on the set $O\in\Bf(\Re)$  if there exist a
measurable function $\varrho$ and  a family of jointly measurable
functions $\{G_{i}:\Re^m\times \Re\to \Re\}$ such that

(i) for every $i$ and $\lambda^m$-almost every $z$, the function
$G_{i}(z,\cdot)$ is an increasing one;

(ii) $\varrho>0$ $\lambda^m$-almost surely and, for every $n\geq
1,z\in\Re^m$, \be\label{9381} |F_n(z+t)-F_n(z+s)|\leq
\sum_{i=1}^m\Big[G_{i}(z, t_i\vee s_i)-G_{i}(z, t_i\wedge
s_i)\Big], \quad \|t\|,\|s\|<\varrho(z), \, z+t\in O,\, z+s\in O.
\ee

It follows from  (\ref{938}), (\ref{934}) and Fubini theorem that,
for every $i=1,\dots,m$ and $P^\Yf$-almost all $\ups\in\Yf$, the
sequence $\{[f^n_i]_\ups\}$ has a uniformly dominated increments
on $[\Omega']_\upsilon$. Since $\lambda^m$-almost every point of
$[\Omega']_\upsilon$ is a density point for $[\Omega']_\upsilon$,
at almost every point of $[\Omega']_\upsilon$ the approximative
partial derivatives of the functions $[f^n_i]_\ups \1_{
[\Omega']_\upsilon}, n\geq 1, i=1,\dots,m$ coincide with those of
the functions $[f^n_i]_\ups, n\geq 1, i=1,\dots,m$.   Thus Theorem
\ref{t218} follows from the Fubini theorem, the statement II of
Proposition \ref{p27} and the following lemma.

\begin{lem}\label{l332} Suppose the sequence $\{F_n:\Re^m\to \Re\}$ to converge $\lambda^m$-a.s.
and to have a uniformly dominated increments on $O\in\Bf(\Re^m)$.
Then the sequence $\{F_n\1_O\}$ is uniformly approximatively
Lipschitz.
\end{lem}

\demo Let $\delta,R>0$ be fixed. Denote $A_{\eps,R}=\{z:\|z\|\leq
R, \varrho(z)> 2m\eps\}$ and take $\eps>0$ such that
$\lambda^m(B_\Re^m(0,R)\backslash A_{\eps,R})<{\delta\over 3}$.
Consider the family of a rectangles of the type
$\prod_{i=1}^m(n_i\eps,(n_i+1)\eps), n_1,\dots, n_m\in \ZZ$. This
family performs a partition of $\Re^m$ up to a set of zero
Lebesgue measure. In this family, consider the sets that provide
non-empty intersections with $A_{\eps,R}$ and denote these sets by
$B^1,\dots, B^J$, here $J<+\infty$ is the total number of the
sets.

Let $j\in\{1,\dots, J\}$ be fixed. Then there exists $z^j\in B^j$
such that $\varrho(z^j)>2m\eps$. Since the diameter of $B^j$ does
not exceed $2m\eps$, this provides that every point  $x\in B^j$
can be written to the form $x=z^j+t$ with $|t|<\varrho(z^j)$.
Denote, for $i=1,\dots,m$, $G^j_i(r)=G_i(z^j, r-z^j_i), r\in\Re$.
Then, from (\ref{9381}), we have that \be\label{9382}
|F_n(x)-F_n(y)|\leq \sum_{i=1}^m [G_i^j(x_i\vee
y_j)-G_i^j(x_i\wedge y_j)],\quad x,y\in B^j\cap O. \ee The set
$B^j$ is a product of intervals $(c_i^j,d_i^j), i=1,\dots, m$.
Every function $G_i^j$ is monotonous and thus differentiable at
$\lambda^1$-almost all points of the interval $(c_i^j,d_i^j)$.
Thus Lemma \ref{l933} provides that, for every $\gamma>0$, there
exists a function $G_{i,\gamma}^j\in C^1(\Re)$ and a compact set
$K_{i,\gamma}^j\subset(c_i^j,d_i^j)$ such that $G_{i}^j$ is
continuous at every point of $K_{i,\gamma}^j$,
$G_{i}^j=G_{i,\gamma}^j$ on $K_{i,\gamma}^j$ and
$\lambda^1((c_i^j,d_i^j)\backslash K_{i,\gamma})^j\leq \gamma$.
Denote $K_{\gamma}^j=\prod_iK_{i,\gamma}^j$. One can choose
$\gamma$ small enough for
$$
\sum_{j=1}^J\lambda^m\left(B^j\backslash
K_{\gamma}^j\right)<{\delta\over 6}.
$$
The function $G_{i,\gamma}^j$ is locally Lipschitz, and therefore
the restriction of $G_{i}^j$ to $K_{i,\gamma}^j$ is Lipschitz with
some constant $L_{i,\gamma}^j$. Then, by (\ref{9382}),  the
restriction of $F_n\1_O$ to $K_{\gamma}^j\cap O$ is Lipschitz with
the constant $L_{\gamma}^j=\sum_{i} L_{i,\gamma}^j$. The
restriction of $F_n\1_O$ to $B(0,R)\backslash O$ is also Lipschitz
with the constant $0$.

By Ulam theorem, there exist compact sets $\hat K\subset \bigcup_j
(K_{\gamma}^j\cap O)$ and $\tilde K\subset B(0,R)\backslash O$
such that
$$
\lambda^m\Big(\bigcup_j (K_{\gamma}^j\cap O)\backslash \hat
K\Big)+ \lambda^m\Big(B(0,R)\backslash (O\cup \tilde
K)\Big)<{\delta\over 6}.
$$
At last, by Egorov theorem, there exist $C>0$ and a compact set
$K_*$ with $\lambda^m(B_\Re^m(0,R)\backslash K_*)< {\delta\over
3}$ such that $|F_n(x)|\leq C, x\in K_*, C>0$. By the
construction, there exists $\theta>0$ such that $\|x-y\|\geq
\theta$ as soon as $x\in K_\gamma^{j_1},y\in K_\gamma^{j_2}$ with
$j_1\not=j_2$ or $x\in \hat K, y\in \tilde K$. Therefore, for
every $n$, the restriction of $F_n\1_O$ to $K_{\delta,R}\eqdef
\hat K\cap\tilde K\cap
K_*\cap\left[\bigcup_{j=1}^JK_{\gamma}^j\right]$ is Lipschitz with
the constant $ L_{\delta,R}\eqdef\max\left[{C\over \theta},\max_j
L_{\gamma}^j\right]. $ By the construction, $A_{\eps,R}\subset
\bigcup_j B^j$ and thus
$$
\lambda^m(B(0,R)\backslash K_{\delta,R})<
\lambda^m(B(0,R)\backslash A_{\eps,R})+{\delta\over
6}+{\delta\over 6}+{\delta\over 3}<\delta.
$$
This completes the proofs of Lemma \ref{l332} and Theorem
\ref{t218}.

\section{Absolute continuity and convergence in variation of distributions to SDE's with jumps}

In this section, applications of Theorems \ref{t212} -- \ref{t218}
to solutions of SDE's with jumps are given. We consider separately
two classes of SDE's. The first one contains \emph{SDE's with
additive noise} of the type (\ref{31}).  The second one contains
\emph{SDE's with non-additive noise}, including SDE's with
non-constant jump rate, of the type (\ref{9311}). The latter class
does not cover the former one because the conditions imposed on
the measure $\nu$ and the coefficients of (\ref{9311}) imply that
the solution to (\ref{9311}) possesses trajectories with bounded
variation, while the L\'evy process $Z$ in (\ref{31}) may be
arbitrary.

Let us introduce notational conventions.  Any time the functional
$f$ of $\nu$ is expressed explicitly through the coefficients
$a,b,c$ and the point measure $\nu$, $f^n$ denotes the functional
of the same form with the coefficients $a^n,b^n,c^n$ and the same
point measure. We introduce conditions $H_1,H_2,\dots$ for a one
functional $f$ in the terms of the coefficients involved into
expression for this functional ($a,b,c$ etc.). Then we write
$H_1^*,H_2^*,\dots$ for the uniform analogues of these conditions,
imposed on the sequence $\{f^n\}$. The constants in these
conditions, as well as the auxiliary functions
$\alpha,\beta,\dots$, are the same with those in conditions
$H_1,H_2,\dots$.   The partial derivative w.r.t. time variable is
denoted by $\prt_t$. The gradient w.r.t. phase variable
$x\in\Re^m$ is denoted by $\nabla$. The unit sphere in $\Re^m$ is
denoted by $S^m$.

\subsection{SDE's with additive noise}

Let $\UU=\UU_1\cup\UU_2$ with $\Pi(\UU_1)<+\infty$. Denote
$$
Z(t)=\int_0^t\int_{\UU_1}c(u)\nu(ds,du)+\int_0^t\int_{\UU_2}c(u)\tilde
\nu(ds,du),\quad t\in\ax,$$ where $c:\UU\to \Re^m,$
$\|c\|\1_{\UU_2}\in L_2(\Pi)$. Consider SDE driven by the {L\'evy
process} $Z$: \be\label{31} X(x,t)=x+\int_0^ta(X(x,s))\, ds+Z(t),
\quad x\in\Re^m, t\in\ax. \ee Under condition

$H_1.$ $a\in C^1(\Re^m,\Re^m), \|a(x)\|\leq C(1+\|x\|),$

\noindent equation (\ref{31}) possesses unique strong solution.
 Put
$f=X(x,t)$, $f^n=X^n(x^n,t^n)$, $\Delta(x,u)=a(x+c(u))-a(x)$,
$$
\Nf(f)=\{\Sf(f)=\Re^m\},\quad
\Sf(f)=\Span\Big\{\Ef_\tau^t\Delta(X(\tau-),p(\tau)),\tau\in\Df\cap[0,t]\Big\},
$$
where $\Ef_s^t, 0\leq s\leq t$ denotes the \emph{stochastic
exponent}, i.e. the $m\times m$-matrix valued process defined by
the equation
$$\Ef_{s}^t=I_{\Re^m}+\int_s^t \nabla
a(X(x,r))\Ef^{r}_{s}\,dr,\quad t\geq s.
$$

\begin{thm}\label{t933} 1. Under condition $H_1$,  $P|_{[X(x,t)]}\circ f^{-1}\ll
\lambda^m.$

2. Let $x^n\to x, t^n\to t$, $c^n(\cdot)\to c(\cdot)$ $\Pi$-almost
everywhere and $a^n\to a, \nabla a^n\to \nabla a$ uniformly on
every compact set. Suppose also that $H_1^*$ holds true and
 $$
 \|c^n(u)\|\1_{\UU_2}\leq \alpha(u),\quad u\in \UU_2\quad \hbox{with}\quad
 \alpha\in L_2(\Pi).
 $$
Then, for every $A\subset{\Nf(f)}$,
$$
P\Bigl|_{A}\circ [X^n(x^n,t^n)]^{-1}\to P\Bigl|_{A}\circ
[X(x,t)]^{-1},\quad n\to +\infty
$$
in variation.

3. Let
 there exist $\eps>0$ such that
\be\label{11}
 \Pi\Bigl(u : (\Delta(y,u),  l)_{\Re^m}\not=0
\Bigr)= +\infty, \quad l\in S^m, y\in \bar B(y,\eps)\eqdef
 \{y:\|y-x\|\leq \eps\}. \ee
Then $P(\Nf(f))=1$.
\end{thm}

\begin{cor}\label{c31} Let condition (\ref{11}) hold true for every $x\in
\Re^m$. Then the transition probability for the process $X$,
considered as a Markov process, possesses a density: $P(X(x,t)\in
dy)=p_{x,t}(y)\, dy$. Moreover, the mapping
$$
\Re^m\times(0,+\infty)\ni (x,t)\mapsto p_{x,t}(\cdot)\in
L_1(\Re^m,\lambda^m)
$$
is continuous and, consequently, $X$ is a strongly Feller process.
\end{cor}

\begin{rem}\label{r41} Examples are available (see \cite{Me_jumps_UMZH},
Example 1.4 and Proposition 1.2), such that $p_{x,t}\not\in
L_{p,loc}(\Re^m)$ for every $p>1, x\in \Re^m, t>0$. This means
that, in some sense, the continuity property exposed in the
Corollary \ref{c31} is the best possible one when no additional
restrictions on the measure $\Pi$ are imposed.
\end{rem}

The proof of Theorem \ref{t933} contains several steps. First, we
prove differentiability of $f$ and the property of $\{f^n\}$ to
have a uniformly dominated increments. Let a grid $\Gf$ of
dimension $m$ be fixed.

\begin{prop}\label{p33} Under conditions of Theorem \ref{t933},  every component of the
vector $X(x,t)$ is a.s. differentiable w.r.t. $\Gf$ for every
$x\in\Re^m,t\in\ax.$ For every $i=1,\dots,m$, the process
$Y_i(x,\cdot)=(D_i^\Gf X_j(x,t))_{j=1}^m$ satisfies the equation
\be\label{32}
 Y_i(x,t)=\int_0^t\int_{\Gamma_i} \Delta(X(x,s-),u
 )Jh_i(s)\,\nu(ds,du)+
\int_0^t [\nabla a](X(x,s)) Y_i(x,s)\,ds ,\quad t\geq 0. \ee
\end{prop}

\begin{rem}\label{r53} In the case $m=1$, the analogous result
  was proved in \cite{Nou_sim}. We cannot use here
 the result from  \cite{Nou_sim} straightforwardly, since the proof
 there contains some specifically one-dimensional features such
 as an exponential formula for the derivative of the flow
 corresponding  to ODE (Lemma 1 \cite{Nou_sim}).
\end{rem}

\demo We fix $i\in \{1,\dots, m\}$ and omit the subscript $i$ in
the notation. Denote $\nu^\Gamma(t,A)=\nu(t,A\backslash \Gamma),$
$$Z^\Gamma(t)=\int_0^t\int_{\UU_1\backslash \Gamma}c(u)
\nu(ds,du)+\int_0^t\int_{\UU_2\backslash \Gamma}c(u)\tilde
\nu(ds,du)- \int_0^t\int_{\UU_2\cap \Gamma}c(u)\, \Pi(du)\, ds.$$
For a given $t>0, \tau\in(0,t),p\in\UU, x\in\Re^m$, consider the
process $X_\cdot^{\tau}$ on $[0,t]$ defined by
$$
X_{r}^{\tau}=\begin{cases}x+\int_{0}^r
a(X^{\tau}_s)\,ds+Z^\Gamma(r),&r<\tau\\
x+\int_{0}^r a(X^{\tau}_s)\,ds+c(p)+ Z^\Gamma(r),&r\geq
\tau\end{cases}.
$$

Denote $\Omega_k=\{\Df\cap\{0,t\}=\emptyset,\#(\Df^\Gamma\cap
(0,t))=k, \},k\geq 0.$ Since $\Gamma\in\pf$,
$\Omega=\bigcup_k\Omega_k$ almost surely. Thus, it is enough to
prove a.s. differentiability on every $\Omega_k$, separately. The
case $k=0$ is trivial, let us consider the case $k=1$.

By the construction of the transformations $T_{th}^\Gamma$,
\be\label{221} [T_{th}^\Gamma
\tau_j^\Gamma](\omega)=T_{-th}(\tau_j^\Gamma(\omega)),\quad
\omega\in\Omega. \ee The point process $\{p(r), r\in\Df^\Gamma\}$
is independent of $\nu^\Gamma$, and the distribution of the
variable $\tau_1^\Gamma= \min \Df^\Gamma$ is absolutely
continuous. Thus a.s. differentiability of $X(x,t)$ on $\Omega_1$
follows immediately from (\ref{221}) and the following lemma.

\begin{lem}\label{l321} With probability 1, for
$\lambda^1$-almost all $\tau\in (0,t)$,
$$
{d\over d\eps}\Bigl|_{\eps=0}X_{t}^{\tau+\eps}=-
\Delta(X_{\tau-}^{\tau},p)\Ef_{t}
$$
with $\Ef_\cdot$ defined by the equation
$$
\Ef_{r}=I_{\Re^m}+\int_\tau^r \nabla
a(X^{\tau}(s))\Ef_{s}\,ds,\quad r\geq \tau.
$$
\end{lem}

\demo $X_{t}^{\tau}$ is the value at the point $t$ of the solution
to the equation \be\label{33} d\tilde X_r=a(\tilde
X_r)\,dr+dZ^\Gamma(r), \ee with the starting point $\tau$ and the
initial value $X_{\tau}^\tau=X_{\tau-}^{\tau}+c(p).$ Let
 $\eps<0$.  Then $X_{s}^\tau=X_{s}^{\tau+\eps}, s<\tau+\eps.$ Thus
  $X_t^{\tau+\eps}$ is also the value of the solution to the same
 equation with the same starting and terminal points and with  the initial value being equal to
$$
X_{(\tau+\eps)-}^{\tau}+c(p)+\int_{\tau+\eps}^\tau
a(X_{s}^{\tau+\eps})\,ds+ [Z^\Gamma(\tau)-Z^\Gamma ({\tau+\eps})].
$$
Thus the difference $\Phi(\tau,\eps)$ between the initial values
for $X_{t}^{\tau+\eps},X_{t}^{\tau}$ is equal to
 $\int_{\tau+\eps}^\tau[a(X_s^{\tau+\eps})-
a(X_s^{\tau})]\,ds$.

The process  $Z^\Gamma$ has  c\`{a}dl\`{a}g trajectories, and
therefore almost surely the set of discontinuities for its
trajectories is at most countable. In addition, every given point
$s\in \ax$ is a continuity point for the trajectory of $Z^\Gamma$
almost surely. Therefore, there exists a set
$\TT=\TT(\omega)\subset\ax$ of the full Lebesgue measure such that
$$
\delta(t,\gamma)\equiv
\sup_{|s-t|\leq\gamma}[\|Z^\Gamma(s)-Z^\Gamma(t)\|]\to 0,\quad
\gamma\to 0,\quad t\in\TT
$$
and $\tau\in \TT$ a.s. Then $
\|X_{s}^\tau-X_{\tau-}^\tau\|+\|X_{s}^{\tau+\eps}-X_{\tau-}^\tau
-c(p)\|\leq \Cd\{|\eps|+\delta(\tau,|\eps|)\} $ for
$s\in(\tau+\eps,\eps)$. Here and below, $\Cd$ denotes any constant
such that it can be expressed explicitly, but its  exact form is
not needed in a further exposition. Thus, for $\tau\in\TT$,
\be\label{942}
\|\Phi(\tau,\eps)+\eps[a(X_{\tau-}^\tau+c(p))-a(X_{\tau-}^\tau)]\|\leq
\Cd|\eps| \{|\eps|+\delta(\tau,|\eps|)\}. \ee The solution to
(\ref{33}) with the starting point $\tau$ is differentiable w.r.t.
initial value with the derivative being equal $\Ef_\cdot$. This
statement is quite standard and we omit the proof. This together
with (\ref{942}) implies the needed statement.

The case $\eps>0$ is analogous, let us discuss it briefly. Again,
take $\tau\in\TT$ and represent $X_{t}^{\tau}$ as the solution to
(\ref{33})  with the initial value $X_{\tau-}^\tau+p$.
$X_{t}^{\tau+\eps}$ is also the solution to (\ref{33}) but with
the other starting point $\tau+\eps$. The estimates analogous to
ones made before show that, up to the $o(|\eps|)$ terms,
$$
X_{\tau+\eps}^{\tau+\eps}-X_{\tau+\eps}^\tau=\eps \Bigl\{
-a(X_{\tau-}^\tau+p) +a(X_{\tau-}^\tau)\Bigr\},
$$
which implies the statement of the lemma. The lemma is proved.

Now let $k>1$ be fixed. Consider the countable family $\Qf_k$ of
partitions $Q=\{0=q_0<q_1\dots<q_k=t\}$ with
$q_1,\dots,q_{k-1}\in\QQ$ and denote
$$
\Omega_Q=\{\Df\cap\{q_i, i=0,k\}=\emptyset,
\Df^\Gamma\cap(q_{i-1},q_i)=1,i=1,\dots,k\},\quad Q\in\Qf_k.
$$
We have $\Omega_k=\cup_{Q\in\Qf_k}\Omega_Q$. Therefore, it is
enough to verify a.s. differentiability of $X(x,t)$ on $\Omega_Q$
for a given $Q$. The distributions of the variables
$\tau_j^\Gamma, j=1,\dots,k$ are absolutely continuous. Then one
can write the statements analogous to the one of Lemma \ref{l321}
on the intervals $[0,q_1],[q_1,q_2],\dots,[q_{k-1},t]$ and obtain
a.s. differentiability  of $X(x,t)$ from (\ref{221}) and the
theorem on differentiability of the solution to (\ref{31}) w.r.t.
initial value. The proposition is proved.

\begin{prop}\label{p34} Let conditions of Theorem \ref{t933} hold
true. Let, in addition,  the sequence $\{a^n\}$ be uniformly
bounded and $Jh_i(t_n)=0,n\geq 1, i=1,\dots, m$. Then the sequence
$\{X^n(x^n,t^n)\}$ has a uniformly dominated increments w.r.t.
$\Gf$ on $\Omega'=\Omega$.
\end{prop}

\demo Again, we omit $i$ in notation. In the framework of Lemma
\ref{l321}, one has the estimate \be\label{34}
\|X_t^{\tau+\eps}-X_t^{\tau}\|\leq \Cd |\eps|,\quad
\tau,\tau+\eps\in (0,t), \ee valid point-wise. Indeed, both
$X_t^{n,\tau+\eps}$ and $X_t^{n,\tau}$ are the solutions to
(\ref{33}) with the same initial point ($\tau$ for $\eps<0$ and
$\tau+\eps$ for $\eps>0$) and different initial values. The
difference between the initial values are estimated by
$$
\left\|\int_{\tau+\eps}^\tau[a(X_s^{\tau+\eps})-
a(X_s^{\tau})]\,ds\right\|\leq -2\|a\|_\infty \eps\hbox{ for }
\eps<0 \hbox{ and }
\left\|\int^{\tau+\eps}_\tau[a(X_s^{\tau+\eps})-
a(X_s^{\tau})]\,ds\right\|\leq 2\|a\|_\infty \eps\hbox{ for }
\eps>0.
$$
Thus inequality (\ref{34}) follows from  the Gronwall lemma. Using
the described above technique, involving partitions $Q\in\Qf_k$,
and applying the Gronwall lemma once again, we obtain that, almost
surely on the set $\Omega_k$, \be\label{943} \|T_{\eps h}^\Gamma
X(x,t)-X(x,t)\|\leq k\Cd \sup_s |Jh(s)||\eps|. \ee Here we have
used that $Jh(t)=0$ and thus $T_{th}x\in (0,t)$ as soon as $x\in
(0,t)$.

The same estimate holds true for every $n$. Thus every sequence
$\{X_j^n(x^n,t^n)\}, j=1,\dots, m,$ satisfies (\ref{938}) with
$\varrho\equiv +\infty$ and $g_i(t)= t\,\Cd\sup_s |Jh_i(s)|
\int_{\Gamma_i}\nu([0,T]\times \Gamma_i), t\in\Re, i=1,\dots, m,
T=\sup_n t_n$. The proposition is proved.

Now we apply Theorems \ref{t212} and \ref{t218} in order to prove
statements 1 and 2 of Theorem \ref{t933}.

\emph{Proof of statement 1 of Theorem \ref{t933}.}  Consider the
family $\{\UU_N, N\geq 1\}$ of bounded measurable subsets of $\UU$
such that $\UU_N\uparrow \UU, N\to +\infty$.  Denote, by $\Sf^N$,
a linear span of the set of vectors $\{\Ef_{\tau}^t\cdot
\Delta(X(\tau-),p(\tau)),\tau\in\Df^{\UU_N}\cap (0,t)\}$ and put
$\Omega^N=\{\omega: \Sf^N(\omega)=\Re^m\}.$ It is clear that
$\Nf(f)=\bigcup_{N\geq 1}\Omega^N$. Thus, in order to prove
statement 1 of Theorem \ref{t933}, it is enough to prove that
$P|_{\Omega^N}\circ [X(x,t)]^{-1}\ll\lambda^m$ for a given $N$.

Let $N$ be fixed. Denote by $L^M_t$ the set of all vectors
$l=(l_1,\dots, l_m)\subset \NN^m$ with $l_1<l_2\dots<l_m$ and
$Ml_m<t$. Consider the family of differential grids $\{\Gf^{M,l},
l\in L^M_t, M\geq 1\}$ of the form $\Gamma_i^{M, l}=\UU_N,$
$$
  a_i^{M,l}={l_i-1\over M}, b_{i}^{M,
l}={l_i\over M}, \quad h_i^{M, l}(s)=h\left({s-a_i^{M, l}\over
b_i^{M, l}-a_i^{M, l}}\right), \quad s\in\ax,\, i=1,\dots, m,\,
l\in L^M_t,\, M\geq 1,
$$
where $h\in H_0$ is some fixed function such that $Jh>0$ inside
$(0,1)$ and $Jh=0$ outside $(0,1)$.

Our aim is to show that almost surely \be\label{35}
\Omega^N\subset \bigcup_{l,M}\,\{\omega:
\Sigma^{X(x,t),\Gf^{M,l}}(\omega) \hbox{ is non-degenerate }\},
\ee see Theorem \ref{t212} for the notation $\Sigma^{f,\Gf}$.
Theorem \ref{t212} together with (\ref{35}) immediately imply the
needed statement.

Denote  $A_M^{N,t}=\Bigl\{\omega: \Df \cap\left\{{i-1\over
N},i\geq 1\right\}=\emptyset,  \#\{\tau\in \Df^{\UU_N}\cap
(a_i^M,b_i^M)\}\subset\{0,1\}, i=1,\dots, [Mt+1]\Bigr\}$. Since
$\UU_N\in\pf$, one has that  almost surely $
\Omega^N\subset\bigcup_M\, [\Omega^N\cap A_M^{N,t}]. $ On the
other hand, $\Omega^N\cap A_M^{N,t}=\bigcup_{l\in L^M_t}
A_{M,l}^{N,t}$ with
$$
A_{M,l}^{N,t}=A_M^{N,t}\cap\left\{\omega:
\Span\left\{\Ef_{\tau}^t\cdot \Delta(X(\tau-),p(\tau)), \tau\in
\Df^{\UU_N}\cap \left({l_i-1\over M},{l_i\over M}\right),
i=1,\dots, m \right\}=\Re^m\right\}.
$$
Thus, in order to prove (\ref{35}), it is sufficient to show that,
for every $M,l$, the matrix $\Sigma^{X(x,t),\Gf^{M,l}}$ is
non-degenerate on the set $\Omega^N\cap A_{M,l}^{N,t}$. By
(\ref{32}),
$$
D_i^{\Gf^{M,l}}X(x,t)=Jh(\tau_i^{M,l})\Ef_{\tau_i^{M,l}}^t\Delta(X(x,\tau_i^{M,l}-),p(\tau_i^{M,l})),\quad
i=1,\dots,m,
$$
on the set $A_{M,l}^{N,t}$, where $\tau_i^{M,l}$ denotes the
(unique) point from $\Df^N\cap \left({l_i-1\over M},{l_i\over
M}\right)$. By the construction, $Jh(\tau_i^{M,l})>0, i=1,\dots,m$
and the family
$$
\Ef_{\tau_i^{M,l}}^t\Delta(X(x,\tau_i^{M,l}-),p(\tau_i^{M,l})),\quad
i=1,\dots,m,
$$
has the maximal rank on the set $\Omega_N\cap A_{M,l}^{N,t}$. Thus
the matrix $\Sigma^{X(x,t),\Gf^{M,l}}$ is non-degenerate on this
set. This completes the proof of statement 1.

\emph{Proof of statement 2 of Theorem \ref{t933}.} Consider first
the case with $\{a^n\}$ uniformly bounded. The standard limit
theorem for SDE's provide that $X^n(x^n,t^n)\to X(x,t)$ in
probability and, for every grid $\Gf$, $Y_i^n(x^n,t^n)\to
Y_i(x,t)$ in probability. Then Theorem \ref{t218} and Proposition
\ref{p34} imply immediately that, for every $B\in\Ff$,
\be\label{9419} P|_{B\cap A_{M,l}^{N,t}}\circ[f^n]^{-1}\tov
P|_{B\cap A_{M,l}^{N,t}}\circ f^{-1},\quad M,N\geq 1, l\in L_t^M.
\ee We have already proved that
$\Nf(f)=\bigcup_{M,N,l}A_{M,l}^{N,t}$, and thus (\ref{9419})
provides the required statement. The additional limitation on
$\{a^n\}$ to be uniformly bounded can be removed via the following
standard localization procedure. Take, for $R>0$, the function
$a_R$ and the uniformly bounded sequence $\{a^n_R\}$ such that
$a^n_R\to a_R,\nabla a^n_R\to \nabla a_R$ uniformly over every
bounded set and $a_R^n(x)=a^n(x), a_R(x)=a(x), \|x\|\leq R$. Then,
on the set $\{\sup_n\sup_{s\leq t_n}\|X^n(x^n,s)\|\leq R\}$,
solutions to (\ref{31}) with the coefficients $a^n$  coincide with
the solutions to (\ref{31}) with the coefficients $a^n_R$,
respectively. Thus, for every $A\subset \Nf(f)$,
$$
P|_{A\cap\{\sup_n\sup_{s\leq t_n}\|X^n(x^n,s)\|\leq
R\}}\circ[f^n]^{-1}\tov P|_{A\cap \{\sup_n\sup_{s\leq
t_n}\|X^n(x^n,s)\|\leq R\}}\circ f^{-1}.
$$
One can see (the proof is standard and omitted) that
$P(\sup_n\sup_{s\leq t_n}\|X^n(x^n,s)\|\leq R\})\to 1, R\to
+\infty$. This completes the proof of statement 2.

\emph{Proof of statement 3 of Theorem \ref{t933}.} We have that,
under conditions of Theorem \ref{t933}, \be\label{37}
\gamma_N\equiv\inf\limits_{y\in \bar B(x,\eps), v\in
S^m}\Pi\Bigl(u\in\UU_N: (\Delta(y,u),v)_{\Re^m}\not=0\Bigr)\to
+\infty,\quad N\to+\infty.\ee This statement follows immediately
from the Dini theorem applied to the monotone sequence of lower
semi-continuous functions
$$
\phi_N: \bar B(x,\eps)\times S^m\ni(y,v)\mapsto
\Pi\Bigl(u\in\UU_N: (\Delta(y,u),v)_{\Re^m}\not=0\Bigr).
$$

With probability 1, the matrix $\Ef_0^r$ is  invertible for every
$r$ and the function $r\mapsto  \Ef_0^r$ is continuous (e.g.
\cite{protter}, Chapter 5, \S 10). In addition,
$\Ef_r^t=\Ef_0^t[\Ef_0^r]^{-1},r\in[0,t]$. Therefore,
$$
\Sf(f)=\Re^m\Longleftrightarrow
\Span\Big\{[\Ef_0^{\tau}]^{-1}\cdot
\Delta(X(\tau-),p(\tau)),\tau\in\Df\cap (0,t)\Big\}=\Re^m
$$
almost surely. Denote by $\SS$ the set of all proper subspaces of
$\Re^m$. This set can be parameterized in such a way that it
becomes a Polish space, and, for every family of random vectors
$\xi_1,\dots,\xi_k$, the map $\omega\mapsto
\Span(\xi_1(\omega),\dots,\xi_k(\omega))$ defines the random
element in $\SS$.

For every $N\geq 1$, consider the set
$\Df^{\UU_N}=\{\tau_1^N,\tau_2^N,\dots\}$. Denote
$$\Sf_\delta^N=\Span\Big\{[\Ef_0^{\tau}]^{-1}\cdot
\Delta(X(\tau-),p(\tau)),\tau\in\Df^{\UU_N}\cap (0,t)\Big\},\quad
\Sf_\delta=\Span(\bigcup_N\Sf_\delta^N).
$$
For a given $\Sf^*\in \SS$, $\delta>0$,  consider the event
$$
D_\delta^N=\{\Sf_\delta^N\not \subset \Sf^*\}=\{\exists k:
\tau_k^N\leq \delta,
[\Ef_0^{\tau_k^N}]^{-1}\Delta(X(\tau_k^N-),p(\tau_k^N))\not\in
\Sf^*\}
$$
One has that $\Omega\backslash D_\delta^N\subset B_\delta\cup
C_\delta^N,$ where $ B_\delta=\{\exists s\in[0,\delta]:
X(s-)\not\in \bar B(x,\eps)\},$
 $$C_\delta^N=\bigcap_k\Bigl[\{\tau_k^N>\delta\}\cup \{X(\tau_k^N-)\in \bar
B(x,\eps),
[\Ef_0^{\tau_k^N}]^{-1}\Delta(X(\tau_k^n-),p(\tau_k^n))\in
\Sf^*,\tau_k^n\leq \delta\}\Bigr].
$$
The distribution of the value $p(\tau_k^n)$ is equal to
$\lambda_N^{-1} \Pi|_{\UU_N}$, where $\lambda_N=\Pi(\UU_N)$.
Moreover, this value is independent with the $\sigma$-algebra
$\Ff_{\tau_k^N-}$, and, in particular, with  $X(\tau_k^N-),
\Ef_0^{\tau_k^N}$. This provides the estimate \be\label{38}
P\Bigl[\{\tau_k^N>\delta\}\cup \{X(\tau_k^N-)\in \bar B(x,\eps),
[\Ef_0^{\tau_k^N}]^{-1}\Delta(X(\tau_k^N-),p(\tau_k^N))\in
\Sf^*,\tau_k^N\leq \delta\}\Bigr|\Ff_{\tau_k^N-}\Bigr]\leq
$$
$$
\leq \1_{\{\tau_k^N>\delta\}}+(1-{\gamma_N\over
\lambda_N})\1_{\{\tau_k^N\leq\delta\}} \ee with $\gamma_N$ defined
in (\ref{37}). It follows from (\ref{38}) that
$$
P(C_\delta^N)\leq E\left(1-{\gamma_N\over
\lambda_N}\right)^{\nu([0,\delta]\times\UU_n)}=\exp\{-\delta\gamma_N\}\to
0,\quad N\to+\infty.
$$
 Then
$D_\delta^N\subset\{\Sf_\delta\not \subset \Sf^*\}$, and almost
surely \be\label{39} \{S_\delta \subset S^*\}\subset B_\delta. \ee
Now we take $\delta<{t\over m}$ and iterate (\ref{39}) on the time
intervals $[0,\delta],[\delta,2\delta],\dots,
[(m-1)\delta,m\delta]$ with $\Sf^*_1=\{0\},
\Sf^*_2=\Sf_\delta,\dots, \Sf_m^*=\Sf_{(m-1)\delta}$ (we can do
this due to the Markov property of $X$). We obtain that
$$
\{\dim \Sf_t<m\}\subset \bigcup_{k=1}^m\{\dim
\Sf_{(k-1)\delta}=\dim \Sf_{k\delta}<m\}\subset B_{m\delta}.
$$
Since $P(B_{m\delta})\to 0, \delta\to 0+$, this provides  that
$P\{\dim \Sf_t<m\}=0$ and completes the proof of Theorem
\ref{t933}.

Condition (\ref{11}) involves both the L\'evy measure of the noise
and the coefficient $a$. In some cases, it would be convenient to
have a more explicit sufficient conditions for (\ref{11}), with
the restrictions on $a$ and $\Pi$ separated one from another.

The first condition is given in the case  $m=1$. Denote
$N(a,z)=\{y\in \Re: a(y)=z\}$.

\begin{prop}\label{p36} Suppose that $\Pi(u:c(u)\not=0)=+\infty$ and
there exists some $\delta>0$ such that
$$
\forall z\in\Re\quad \#\Bigl[N(a,z)\cap
(x-\delta,x+\delta)\Bigr]<+\infty.
$$
Then (\ref{11}) holds true, and therefore $P(\Nf(f))=1$.
\end{prop}

\begin{rem} In \cite{Nou_sim},  the law of the solution to
one-dimensional SDE (\ref{31}) was proved to be absolutely
continuous under condition that $a(\cdot)$ is strictly monotonous
at some neighborhood of $x$. One can see that this condition is
somewhat more restrictive than the one of Proposition \ref{p36}.
\end{rem}

{\it Proof of the Proposition.} Take $\eps={\delta\over 2}$. Then,
for every $y\in \bar B(x,\eps)$,
$$\{u:\Delta(y,u)=0\}=\{u: a(y+c(u))=a(y)\}\subset
$$
$$\subset  \{u:|c(u)|>\Cd \delta\}\cup\{u:x+c(u)\in N(a,a(y))\cap
(x-\delta,x+\delta)\}=\Delta_1\cup\Delta_2.
$$
Here we have used that $a$ is Lipschitz. We have that, for every
$d>0$, the set $\{u:|c(u)|>d\}$ has finite measure $\Pi$, and
therefore $\Pi(\Delta_1)<+\infty$. The set $N(a,a(y))\cap
(x-\delta,x+\delta)\backslash \{x\}$ is finite and, therefore,
separated from $x$. Thus
$\Pi(\Delta_2\backslash\{u:c(u)=0\})<+\infty$. Since
 $\Pi(u:c(u)\not=0)=+\infty$,
this means that $\Pi(\Delta(x,u)\not=0)=+\infty$. The proposition
is proved.

The second sufficient condition is formulated in the
multidimensional case. Define a \emph{proper smooth surface}
$S\subset\Re^m$ as any set of the type $S=\{x:\phi(x)\in L\}$,
where $L$ is a proper linear subspace of $\Re^m$ and $\phi\in
C^1(\Re^m,\Re^m)$ is such that $\det \nabla \phi(0)\not=0$ and
$\phi^{-1}(\{0\})=\{0\}$.

\begin{prop}\label{p37} Suppose that one of the following group of conditions holds true:

{$\mathbf{a.}$} $a\in C^1(\Re^m, \Re^m)$, $\det \nabla
a(x_*)\not=0$ and \be\label{310} \Pi(u:c(u)\in \Re^m\backslash
S)=+\infty\quad  \hbox{for every proper smooth surface $S$;} \ee

{$\mathbf{b.}$} $a(x)=Ax, A\in \Lf(\Re^m,\Re^m)$ is non-degenerate
and \be\label{311} \Pi(u:c(u)\in \Re^m\backslash L)=+\infty\quad
\hbox{for every proper linear subspace $L\subset \Re^m$.}\ee

Then (\ref{11}) holds true, and therefore $P(\Nf(f))=1$.
\end{prop}

{\it Proof.} Consider the set $\Phi_{x,\eps}$ of the functions
$\phi_y:\Re^m\ni h\mapsto a(y+h)-a(y)\in \Re^m, y\in \bar
B(x,\eps)$. It is easy to see that if
$$
\Pi(u: c(u)\not\in \phi^{-1}(L_v))=+\infty,\quad \phi\in
\Phi_{x,\eps}\quad  \hbox{for every linear subspace $L_l\equiv
\{y:(y,l)=0\}, l\in S^m$,}
$$
then (\ref{11}) holds true. In the case {\bf b}, $\Phi_{x,\eps}$
contains the unique function $\phi(h)=Ah$. Since $A$ is
non-degenerate, $\phi^{-1} (L_l)$ is a proper linear subspace of
$\Re^m$ for every $v\in S^m$, and (\ref{311}) provides (\ref{11}).
In the case {\bf a}, $\phi_y\in C^1(\Re^m, \Re^m)$, and for $\eps$
small enough $\det \nabla \phi_y(0)=\det \nabla a(y)\not=0, y\in
\hat B(x,\eps)$. Then $\phi^{-1}_y (L_l)$ is a proper smooth
surface for every $l\in S^m$, and (\ref{310}) provides (\ref{11}).
The proposition is proved.

Condition (\ref{310}) holds true, for instance, if $\Pi(u:c(u)\in
\Re^m\backslash Y)=+\infty$ for every set $Y\subset \Re^m$, whose
Hausdorff dimension does not exceed $m-1$. Condition (\ref{311})
is close to the necessary one, this is illustrated by the
following simple example. Let (\ref{311}) fail for some $L$, and
let $L$ be invariant for $A$. Then, for $x\in L$ and any $t\geq
0$, $P(X(x,t)\in L)>0.$ Therefore, the law of $X(x,t)$ is not
absolutely continuous.

Condition (\ref{311}) was introduced by M.Yamazato in
\cite{yamazato}, where the problem of the absolute continuity of
the distribution of the L\'evy process was studied. This condition
obviously is  necessary for the law of $Z(t)$ to possess a
density. In \cite{yamazato}, some sufficient conditions were also
given. Statement 4 of the main theorem in \cite{yamazato}
guarantees the absolute continuity of the law of $Z(t)$ under the
following  three assumptions:

(a) condition (\ref{311}) is valid;

(b) $\Pi(u: c(u)\in L)=0$ for every linear subspace
$L\subset\Re^m$ with dim\,$L\leq m-2$;

(c)  the conditional distribution of the radial part of some
\emph{generalized polar coordinate} is absolutely continuous.

We remark that assumption (c) is some kind of a "spatial
regularity" assumption and is crucial in the framework of
\cite{yamazato}. Without such an assumption, condition (\ref{311})
is not strong enough to guarantee $Z(t)$ to possess a density,
this is illustrated by the following example.

\begin{ex}\label{e39} Let $\UU=\Re^2\backslash \{0\}, c(u)=u,
m=2, \Pi=\sum_{k\geq 1}\delta_{z_k}$, where $z_k=({1\over k!},
{1\over (k!)^2}), k\geq 1.$ Every point $z_k$ belongs to the
parabola $\{z=(x,y): y=x^2\}$. Since every line intersects this
parabola at most at  two points, condition (\ref{311}) and
assumption (b) given before hold true. On the other hand, for any
$t>0$, it is easy to calculate the Fourier transform of the first
coordinate $Z_1(t)$ of $Z(t)=(Z_1(t),Z_2(t))$ and show that
$$
\lim_{N\to +\infty} E\exp\{i 2\pi N! Z_1(t)\}=1.
$$
 This means that the law of $Z(t)$ is singular.
\end{ex}

Although condition (\ref{311}) is not strong enough to provide the
Levy process $Z$ itself to possess an absolutely continuous
distribution, Proposition \ref{p37} shows that this condition
appears to be a proper one for the solution to an
Orstein-Uhlenbeck type SDE driven by this process to possess a
density as soon as the drift coefficient is non-degenerate. At
this time, we cannot answer the question whether (\ref{311})  is
strong enough to handle the non-linear case, i.e. whether
statement {\bf a} of Proposition \ref{p37} is valid with
(\ref{310}) replaced by (\ref{311}).

\subsection{Solutions to SDE's with non-additive noise and non-constant jump rate} Suppose $\UU$
to have the form $\UU=\VV\times\ax$ and the measure $\Pi$ to have
the  form $\Pi=\pi\times \lambda^1$.  Denote $\nu(dt,du)\eqdef
\nu(dt,dv,dp), u\cong(v,p)$ and consider SDE of the type
\be\label{9311} X(x,t)=x+\int_0^ta(X(x,s))\,
ds+\int_0^t\int_\VV\int_0^{b(X(x,s-),v)}c(X(x,s-),v)\nu(ds,dv,dp),\quad
x\in\Re^m, t\in\ax. \ee

 The following conditions are imposed.

$H_2$. $a\in C^1_b(\Re^m,\Re^m)$, $b(\cdot,v)\in C_b(\Re^m,
\Re^+)$, $c(\cdot,v)\in C^1_b(\Re^m, \Re^m)$, $v\in \VV$.

$H_3$. There exist  $\beta,\gamma:\VV\to \Re^+$ such that
$\beta\gamma\in L_1(\VV,\pi)$ and
$$b(x,v)\leq \beta(v),\quad \|c(x,v)\|\leq
\gamma(v),\quad x\in\Re^m, v\in \VV,$$
$$|b(x,v)-b(y,v)|\leq \|x-y\| \beta(v), \quad \|c(x,v)-c(y,v)\|\leq
\|x-y\| \gamma(v), \quad x,y\in\Re^m, v\in \VV.
$$

$H_{4}$. There exists a representation $b(x,v)=b_0(v)+b_1(x,v)$
such that $b_0\gamma\in L_1(\VV,\pi)$ and, for some
$\beta_1,\gamma_1:\VV\to \Re^+$,  ${\beta\beta_1\gamma_1}\in
L_2(\VV,\pi)$  and
$$
|b_1(x,v)|\leq \beta_1(v), \quad \|c(x,v)\|\leq \gamma_1(v), \quad
x,y\in\Re^m, v\in \VV.
$$

$H_5$. $a^n\to a, \nabla a^n\to \nabla a$ and, for $\pi$-almost
all $v\in \VV$, $b^n(\cdot,v)\to b(\cdot,v),\, c^n(\cdot,v)\to
c(\cdot,v), \nabla c^n(\cdot,v)\to \nabla c(\cdot,v)$ uniformly on
every compact set.

Under  conditions $H_2,H_3$, equation (\ref{9311})  possesses
unique strong solutions being a strong Markov processes with
c\'adl\'ag trajectories. Moreover, under conditions $H_2^*,H_3^*,
H_5$, $X^n(x^n, t^n)\to X(x,t)$ in probability for any sequences
$x_n\to x$ and $t_n\to t$ (recall that $X^n$ denotes the solution
to (\ref{9311}) with the coefficients $a,b,c$ replaced by
$a^n,b^n,c^n$). We omit the proofs of these statements, referring
to \cite{fourn_2002}, Section 2 for the proof of a similar
statement.

Put $f=X(x,t)$. Denote by $p_1(\cdot), p_2(\cdot)$ the projections
of the point process $p(\cdot)$ on the first and second
coordinates in $\UU=\VV\times \Re^+$, respectively. For $y\in
\Re^m$, denote $\VV_y=\{v\in \VV: I_{\Re^m}+\nabla c(y,v)$ is
invertible\}, $\Pi_y(dv)=b(y,v)\pi(dv)$ and put
$$
\Delta(y,v)=a(y+c(y,v))-a(y)-\nabla c(y,v)a(y,v), \quad \tilde
\Delta(y,v)=[I_{\Re^m}+\nabla c(y,v)]^{-1}\Delta(y,v),\quad v\in
\VV_y,
$$
$$
\Nf(f)=\{\Sf(f)=\Re^m\}, \quad
\Sf(f)=\Span\left\{\Ef_\tau^t\Delta(X(\tau-),p_1(\tau)),\tau\in\Df\cap[0,t]:
p_2(\tau)\in\Big[0,b(X(\tau-), p_1(\tau))\Big]\right\},
$$
where the stochastic exponent $\Ef_s^t, 0\leq s\leq t$ is defined
by the equation \be\label{9321}\Ef_{s}^t=I_{\Re^m}+\int_s^t \nabla
a(X(x,r))\Ef^{r}_{s}\,dr+
\int_s^t\int_\VV\int_0^{b(X(x,r-),v)}\nabla
c(X(x,r-),v)\Ef^{r-}_{s}\nu(ds,dv,dp),\quad t\geq s. \ee
\begin{thm}\label{t934} 1. Under conditions $H_2,H_3,$
$$P|_{[X(x,t)]}\circ f^{-1}\ll \lambda^m.$$

2. Under conditions $H_2^*$ -- $H_4^*$, $H_5$, for any sequences
$x^n\to x, t^n\to t$, $$P\Bigl|_{A}\circ [X^n(x^n,t^n)]^{-1}\to
P\Bigl|_{A}\circ [X(x,t)]^{-1}$$ in variation for every
$A\subset{\Nf(f)}$.

3. Suppose that, for every $y\in\Re^m,l\in S^m$,
 \be\label{13}
 \Pi_y\Bigl(v\in \VV_y :  (\tilde \Delta(y,v),  l)_{\Re^m}\not=0
\Bigr)= +\infty. \ee Then $P(\Nf(f))=1$.
\end{thm}

\begin{rem} The statements 2 of Theorems \ref{t933},\ref{t934} can be used
efficiently in order to provide the \emph{local Doeblin condition}
to hold true for the Markov processes $X$, see \cite{Me_ergodic}.
In such a set up, the sufficient conditions for $P(\Nf(f))>0$ are
required rather than the conditions for $P(\Nf(f))=1.$ Here we
formulate one condition of such a type: \be\label{12}
 \Pi_x\Bigl(v\in\VV_x : (\tilde \Delta(x,u),  l)_{\Re^m}\not=0, \|c(x,v)\|<\eps
\Bigr)>0, \quad l\in S^m, \eps>0. \ee We do not give the proof
here, referring to the similar proof of Proposition 4.3
\cite{Me_ergodic}. See also Proposition 4.8 \cite{Me_ergodic} for
a refinement of condition (\ref{12}) in the one-dimensional case.
\end{rem}

\begin{rem} Theorem \ref{t934} still holds true with
the uniform bounds on $a,b,c,\nabla c$ replaced by the linear
growth conditions
$$
\|a(x)\|\leq L(1+\|x\|),\quad b(x,v)\leq (1+\|x\|)\beta(v),\quad
\|c(x,v)\|\leq (1+\|x\|)\gamma(v),
$$
$$
b(x,v)\leq (1+\|x\|)\beta_1(v), \quad \|c(x,v)\|\leq
(1+\|x\|)\gamma_1(v).
$$
One can prove this via the localization procedure analogous to the
one used in the proof of statement 2 of Theorem \ref{t933}.
\end{rem}

Conditions $H_3, H_4$  cover a large variety of SDE's, let us
emphasize  some particular classes of equations.

\textbf{A.} Let $b(x,v)=b_0(v),$ then (\ref{9311}) is an SDE with
constant jump rate. The L\'evy measure of the noise is given by
$\Pi'(dv)=b(v)\pi(dv)$ and condition $H_4$ holds true with
$\beta_1\equiv 0, \gamma_1=\gamma$. Condition $H_3$ now has the
form
$$ \|c(x,v)\|\leq
\gamma(v),\quad \quad \|c(x,v)-c(y,v)\|\leq \|x-y\| \gamma(v),
\quad x,y\in\Re^m, v\in \VV,\quad \gamma\in L_1(\VV,d\Pi').
$$
For such an SDE's, Theorem \ref{t934} is already proved in
\cite{Me_TViMc} (statements 1 and 3) and \cite{Me_conv_var}
(statement 2). Also, in this case statement 1 is closely related
to  Theorem 3.3.2 \cite{Denis}, but the later theorem has the
"gap" in its proof, discussed in subsection 4.3 below.

\textbf{B.} Let $\sup_x|b_1(x,v)|\in L_1(\VV,\pi),$ i.e. the jump
rate varies moderately, in a sense.   In this case, $H_3$ yields
$H_4$ with $\beta_1(v)=\sup_x|b_1(x,v)|$ and $\gamma_1=\gamma$. A
class of equations satisfying, among others, the condition
analogous to the one indicated above is studied in \cite{bally}
(the so called  \emph{case without blow up}).

\textbf{C.} Let $b(\cdot,v)=b(\cdot)\in C_b^1(\Re^m)$.   Such
class of (one-dimensional) equations is studied in
\cite{fourn_2008}. In this case,  one can put
$\beta=\beta_1=\,$const, $\gamma=\gamma_1$ and claim $\gamma\in
L_1(\VV,\pi)\cap L_2(\VV,\pi)$.

We remark that in \cite{bally} and \cite{fourn_2008}, for the
cases \textbf{B} and \textbf{C} respectively, existence of a
\emph{smooth} distribution density for the solution to
(\ref{9311}) is proved (see also references therein for some
previous results on absolute continuity of the law of the
solution). This is an essentially stronger result than statement 1
of Theorem \ref{t934}, but the conditions, imposed in \cite{bally}
and \cite{fourn_2008}, are much more restrictive. This is
substantial, because the solution to (\ref{9311}) may possess a
distribution density, but this density may fail to be smooth (see
Remark \ref{r41} and \cite{Me_TViMc}, Section 5). We turn the
reader's attention to the fact that the convergence in variation
holds true under the same weak assumptions that provide absolute
continuity (statement 2 of Theorem \ref{t934}). This allows one to
study ergodic properties of the solution to (\ref{9311}),
considered as a Markov process, under these weak assumptions
(\cite{Me_ergodic}).

The main difficulty in the proof of Theorem \ref{t934} is to get
the differentiability properties, analogous to those given by
Propositions \ref{p33},\ref{p34}. We expose this step in details
and then sketch the rest of the proof. In order to get an
analogues of Propositions \ref{p33},\ref{p34}, we have to
establish the  properties of the solution to (\ref{9311}),
considered as a function of $x$. Consider SDE of the type
(\ref{9311}) with the starting time moment $s$: \be\label{9415}
X(x,s,t)=x+\int_s^ta(X(x,s,r))\,
ds+\int_s^t\int_\VV\int_0^{b(X(x,s,r-),v)}c(X(x,s,r-),v)\nu(dr,dv,dp),\quad
x\in\Re^m, t\geq s\geq 0. \ee We write $\Ef_\cdot^\cdot$ for the
solution to equation of the type (\ref{9321}) with $X(x,r)$
replaced by $X(x,s,r)$.

\begin{lem}\label{l934} Under conditions $H_2 - H_4$, the following properties hold true.

1. For every $T\in\ax$, there exists $\Cs\in\Re^+$ such that
$$
E\|X(x,s,t)-X(y,s,t)\|\leq \Cs\|x-y\|,\quad x\in\Re^m, \quad 0\leq
s\leq t\leq T.
$$

2. There exists an increasing process $\eta(\cdot)$ such that
$\|X(x,s,t_1)-X(x,s,t_2)\|\leq \eta(t_2)-\eta(t_1), s\leq t_1\leq
t_2$.

3. For every $x\in \Re^m, s\leq t$, the function
$X(\cdot,\cdot,\cdot)$ is differentiable w.r.t. every variable at
the point $(x,s,t)$ with probability 1 and \be\label{947} {\prt
X\over \prt t}(x,s,t)=a(X(x,s,t)),\quad  {\prt X\over \prt
x}(x,s,t)=\Ef_s^t,\quad {\prt X\over \prt s}(x,s,t)=-\Ef_s^ta(x).
\ee
\end{lem}

We remark that the function $X(\cdot,\cdot,\cdot)$ may fail to
possess a continuous trajectories. The situation here is like the
one for the Poisson process $N$: the trajectories $\ax\ni t\mapsto
N(t)$ are a.s. discontinuous, but $N'(t)=0$ a.s. for every fixed
$t\in\ax$.

\demo Denote $A=\sup_x\|a(x)\|+\sup_x\|\nabla a\|$. We have
$$
\|X(x,s,t)-X(y,s,t)\|\leq \|x-y\|+A\int_s^t
\|X(x,s,r)-X(y,s,r)\|\, dr+
$$
$$
+\int_s^t\int_\VV
\int_0^{b(X(x,s,r-),v)}\|X(x,s,r-)-X(y,s,r-)\|\gamma(v)\nu(dr,dv,dp)+\int_s^t\int_\VV
\int_{b(X(x,s,r-),v)}^{b(X(y,s,r-),v)}\gamma(v)\nu(dr,dv,dp),
$$
here we have used the notation $\int_b^a=\int_a^b, a<b$. Put
 $D(x,y,t)\eqdef \sup_{0\leq s\leq r\leq
t}E\|X(x,s,r)-X(y,s,r)\|$ and take the expectation in the previous
inequality. Then we have
$$
D(x,y,t)\leq \|x-y\|+L\int_0^t D(x,y,s)\, ds+2\int_0^t
D(x,y,s)\int_\VV \beta(v)\gamma(v)\pi(dv)ds.
$$
Now the statement 1 follows from the Gronwall lemma.

The statement 2 obviously holds true with $\eta(t)=A
t+\int_0^t\int_\VV\int_0^{\beta(v)}\gamma(v)\nu(dr,dv,dp)$. The
second summand $\eta_2(\cdot)$ in the expression for $\eta(\cdot)$
is a L\'evy process with almost all its trajectories being a
singular functions with locally bounded variation. Then Lebesgue
theorem combined with Fubini theorem provides that, for
$\lambda^1$-almost all $t\in\ax$, $\eta_2'(t)=0$ almost surely.
Since $\eta_2$ is time homogeneous, this yields that
$\eta_2'(t)=0$ almost surely for every $t\in\ax$. This provides
the first relation in (\ref{947}). In order to prove the second
and the third relations in (\ref{947}), we need an auxiliary
construction. In order to shorten the notation, we suppose $s=0$
and omit $s$ in the notation.

Denote $\theta=\sqrt{\beta\beta_1\gamma_1}\in L_1(\VV,\pi)$. For a
given $\eps>0$, consider the random set
$$D_{x,t}^\eps=
\{(r,v,p): r\in[0,t], |p-b(X(x,r-),v)|\leq \eps\theta (v)\}\subset
\ax\times \VV\times \ax$$ and put
$$
\Omega_{x,t}^\eps=\{\omega: \forall \tau\in \Df,
(\tau,p_1(\tau),p_2(\tau))\not\in D_{x,t}^\eps\}.
$$
Consider the sequence $\VV_n\uparrow\VV$ with $\pi(\VV_n)<+\infty$
and denote $\Df^n=\{\tau\in \Gf: p_1(\tau)\in \VV_n, p_2(\tau)\in
[0,n]\}$,
$$
\Omega_{x,t}^{\eps,n}=\{\omega: \forall \tau\in \Df^n,
(\tau,p_1(\tau),p_2(\tau))\not\in D_{x,t}^\eps\}.
$$
With probability 1, the set $\Df^n$ can be represented as
$\Df^n=\{\tau_j^n, j\geq 1\}, \tau_1^n<\tau_2^n<\dots$ and $$
\Omega_{x,t}^{\eps,n}=\bigcap_{j=1}^\infty B_j,\quad
B_j\eqdef\{\tau_j^n>t\}\cup\{\tau_{j}^n\leq t,
|p_2(\tau_j^n)-b(X(x,\tau_j^n-))|> \eps\theta
(p_1(\tau_j^n))\},\quad j\geq 1.$$ Denote, by $\{\Ff_t\}$, the
flow of $\sigma$-algebras generated by $X(x,\cdot)$. For every
$j$, the variable $X(x,\tau_j^n-)$ is $\Ff_{\tau_j^n-}$ --
measurable. On the other hand, the variables
$p_1(\tau_j^n),p_2(\tau_j^n)$ are jointly independent on
$\Ff_{\tau_j^n-}\bigvee \sigma(\tau_j^n)$ and have their
distributions equal ${\pi(\cdot\cap\VV_n)\over \pi(\VV_1)}$ and
${\lambda^1(\cdot\cap[0,n])\over n}$, respectively. For every
given $z\in\Re^+$ and every $v\in\VV_n$, the set $\{p\in [0,n]:
|p-z|\leq \eps\theta(v)\}$ is the interval of the length $\leq
2\eps$. Then the probability for $(p_1(\tau_j^n),p_2(\tau_j^n))$
to satisfy the relation $|p_2(\tau_j^n)-z|\leq \eps\theta
(p_1(\tau_j^n))$ does not exceed ${2\eps\over n\pi(\VV_n)}
\int_{\VV_n}\theta(v)\, dv$, and thus \be\label{944}
P(B_j|\Ff_{\tau_{j}^n-})\geq \1_{\tau_{j}^n>t}+\1_{\tau_{j}^n\leq
t}\left[1-{2\eps\over n\pi(\VV_n)} \int_{\VV_n}\theta(v)\,
dv\right], \quad j\geq 1.\ee
 Since $B_{j-1}$ is $\Ff_{\tau_j^n-}$ --
measurable for $j>1$, (\ref{944}) imply the estimate
$$
P(\Omega_{x,t}^{\eps,n})\geq E\left[1-{2\eps\over n\pi(\VV_n)}
\int_{\VV_n}\theta(v)\, dv\right]^{\nu([0,t]\times \VV_n\times
[0,n])}=\exp\left[-2t\eps \int_{\VV_n}\theta(v)\, dv\right].
$$
After passing to the limit as $n\to +\infty$, we get
\be\label{945} P(\Omega_{x,t}^\eps)\geq \exp\left[-2t\eps
\int_{\VV}\theta(v)\, dv\right].\ee Therefore, for every
$x\in\Re^m$ and  $t\in\ax$, almost every $\omega$ belongs to some
$\Omega_{x,t}^\eps$ with $\eps>0$.

Consider the linear SDE $E(r)= 1+\int_0^rE(s-)d\eta(s)$ and write
$\Ls=\Ls(t,\omega)=E(t,\omega)<+\infty$ a.s. Write
$$
\zeta(\kappa,t)=\int_0^t\int_{\theta(v)\leq
\kappa\beta(v)}\int_{(b_0(v)-\beta_1(v))\vee
0}^{b_0(v)+\beta_1(v)}{\gamma_1(v)\beta(v)\over
\theta(v)}\nu(ds,dv,dp),\quad \kappa,t\in\ax.
$$
We have
$$
E\zeta(\kappa,t)\leq 2t \int_{\theta(v)\leq
\kappa\beta(v)}{\beta(v)\beta_1(v)\gamma_1(v)\over
\theta(v)}\pi(dv)=2t \int_{\theta(v)\leq
\kappa\beta(v)}\theta(v)\pi(dv)\to 0,\quad \kappa\to 0.
$$
Since $\zeta(\cdot,t)$ is monotonous, this provides that
$\zeta(\kappa,t)\to 0,\kappa\to 0+$ almost surely.

 Let $x\in\Re^m, t\in\ax,
\eps>0$ and $\omega\in \Omega_{x,t}^\eps$ be fixed. Suppose that
$\Ls(t,\omega)<+\infty$ and $\zeta(\kappa,t,\omega)\to 0,
\kappa\to 0+$. For a given $y\in\Re^m$, we consider  $X(y,\cdot)$
and put $\varsigma=\inf\{r: \|X(x,r)-X(y,r)\|>{2\Ls\|x-y\|}\}$. We
remark that the random variable $\varsigma$ is not a stopping time
since $\Ls$ is defined through the whole configuration of $\nu$.
All the integrals over $\nu$  throughout the rest of the proof
should  be understood, for the fixed $\omega$, in the point-wise
sense. We have \be\label{948}
X(y,r)-X(x,r)=(y-x)+\int_0^r\Big(a(X(y,s))-a(X(y,s))\Big) ds+
$$
$$
+\int_0^r\int_\VV\int_0^{b(X(x,s-),v)}\Big(c(X(y,s-),v)-c(X(x,s-),v)\Big)\nu(ds,dv,dp)+
$$
$$
+\int_{0}^r\int_{\VV}\int_{b(X(x,s-),v)}^{b(X(y,s-),v)}c(X(y,s-),v)\nu(ds,dv,dp).
\ee The last integral in (\ref{948}), for every $\delta>0$, is
dominated by $I_1(\delta,r)+I_2(\delta,r)$,
$$
I_1(\delta,r)=\int_{0}^r\int_{\eps\theta(v)\leq \Ls
\delta\beta(v)}\int_{(b_0(v)-\beta_1(v))\vee
0}^{b_0(v)+\beta_1(v)}\gamma_1(v)\nu(ds,dv,dp),
$$
$$I_2(\delta,r)=\int_{0}^r\int_{\eps\theta(v)> \Ls
\delta\beta(v)}\int_{b(X(x,s-),v)}^{b(X(y,s-),v)}\gamma_1(v)\nu(ds,dv,dp).$$
Take $\delta=2\|y-x\|$. Then, for $s\leq \varsigma$, we have
$\|X(x,s-)-X(y,s-)\|\leq \Ls\delta$. Thus, as soon as
$\eps\theta(v)> \Ls \delta\beta(v)$,
$$
|b(X(x,s-),v)-b(X(y,s-),v)|\leq \Ls \delta\beta(v)<\eps\theta(v).
$$
This means that $I_2(\delta,r)=0, r\leq \varsigma$ on the set
$\Omega_{x,t}^\eps$.

For $\eps\theta(v)\leq \Ls\delta\beta(v)$, we have
$\delta\cdot{\Ls\beta(v)\over\eps\theta(v)} \geq 1$, and
therefore, for $r\leq t$,  $I_1(\delta,r)$ can be estimated by
${\Ls\delta\over\eps}\zeta({\Ls\delta\over\eps}, t).$ Since
$\zeta(\kappa,t)\to 0,\kappa\to 0+$, there exists $\delta_0>0$
such that ${\Ls\over\eps}\zeta({\Ls\delta\over\eps}, t)< {1\over
3}, \delta<\delta_0$.

The first two integrals in (\ref{948}) are dominated by
$$
\int_0^r\|X(x,s-)-X(y,s-)\|d\eta(s).
$$
 Therefore, if $2\|y-x\|<\delta_0$, then,  up to
the time moment $\varsigma$, the values of the process
$\|X(x,\cdot)-X(y,\cdot)\|$ are dominated by the solution to the
equation
$$
Z(r)={5\over 3}\|x-y\|+\int_0^rZ(s-)d\eta(s).
$$
But $Z(r)={5\over 3}\|x-y\|E(r)< 2\Ls\|y-x\|, r\leq t$. This means
that $\varsigma>t,$ and the estimates given before show that
\be\label{9420} \Delta(x,y,t)\eqdef\sup_{r\leq
t}\Big[I_1(\|y-x\|,r)+I_2(\|y-x\|,r)\Big]=o(\|y-x\|),\quad y\to x.
\ee Now we fix $X(x,\cdot)$ and consider (\ref{948}) as a family
of the equations on $X(y,\cdot)$ with the parameter $y$ involved
both into the initial value and the small interfering term (i.e.,
the last summand in (\ref{948}))  dominated by $\Delta(x,y,t)$.
Since the coefficients of these equations are smooth, (\ref{9420})
and  the standard considerations provide that the solution is
differentiable w.r.t. $y$ at the point $y=x$ and the derivative is
equal $\Ef_0^t$.

At last, the discussed above differentiability properties of the
process $\eta$
 provide that, for every $x\in\Re^m, s>0$,
$$
{1\over\Delta s}\Big(X(x,s,s+\Delta s) - x\Big)\to a(x),\quad
{1\over\Delta s}\Big(X(x,s-\Delta s,s) - x\Big)\to a(x),\quad
\Delta s\to 0+
$$
with probability 1. This, together with the already proved
relation ${\prt X\over \prt x}(x,s,t)=\Ef_s^t$, provides that
${\prt X\over \prt s}(x,s,t)=-\Ef_s^ta(x)$ with probability 1. The
lemma is proved.

Inequality $\varsigma>t$ can be rewritten to the following form:
$$\|X(x,r)-X(y, r)\|\leq 2\Ls \|x-y\|,\quad \|x-y\|\leq {1\over
2}\delta_0,\quad r\in[0,t]. $$ We remark that the same arguments
with  those used in the proof of $\varsigma>t$ can be applied, for
the same $\omega, \Ls, \delta_0$, to SDE of the type (\ref{9415})
with the initial value $x'=X(x,0,s)$. These estimates provide
that, for every $s\in[0,t]$,
 \be\label{949}
\sup_{r\in[s,t]}\|X(x,0,t)-X(y,s,t)\|\leq 2\Ls\|X(x,
0,s)-y\|\quad\hbox{as soon as}\quad  \|X(x, 0,s)-y\| \leq {1\over
2}\delta_0. \ee In order study the properties of the sequence of
the solutions to SDE's of the type (\ref{9311}), we need the
following uniform version of (\ref{949}).

\begin{lem}\label{l935} Let conditions $H_2^*$ -- $H_4^*, H_5$
hold true and $x^n\to x, t^n\to t$ be arbitrary sequences. Then
there exist an a.s. positive random variable $\sigma$, a random
variable $\zeta$ and a subsequence $\{n(k), k\geq 1\}\subset \NN$
such  that, for every $k\geq 1, s\in[0, t^{n(k)}],$
\be\label{9410} \sup_{r\in[s,t^{n(k)}]}\|X^{n(k)}(x^{n(k)},
0,r)-X^{n(k)}(y,s,r)\|\leq \zeta\|X^{n(k)}(x^{n(k)},
0,s)-y\|,\quad \|X^{n(k)}(x^{n(k)}, 0,s)-y\| \leq \sigma. \ee
\end{lem}
\demo In order to shorten exposition, we consider the case
$t^n\equiv t$; one can easily see  that this restriction is not
essential in the considerations given below. We omit the starting
moment in the notation and write $X(x,t)$ for $X(x,0,t)$. First,
let us show briefly that \be\label{9411}\sup_{r\leq
t}\|X^n(x^n,r)-X(x,r)\|\to 0\quad \hbox{in probability.}\ee
Estimates analogous to those given at the beginning of the proof
of Lemma \ref{l934} provide that $\sup_{r\leq
t}E\|X^n(x^n,r)-X(x,r)\|\to 0$. We write $$
X^n(x^n,r)=x^n+\int_0^r\tilde a^n(X^n(x^n,s))\, dr+M^n(r),\quad
M^n(r)\eqdef\int_0^r\int_\VV\int_0^{b^n(X^n(x^n,s-),v)}\!\!\!\!c^n(X^n(x^n,s-),v)\tilde\nu(ds,dv,dp),
$$
$$\tilde a^n\eqdef a^n+\int_\VV\int_0^{b^n(\cdot,v)}c^n(\cdot,v)dp\pi(dv)=
a^n+\int_\VV b^n(\cdot,v)c^n(\cdot,v)\pi(dv).
$$
Let $M,\tilde a$ be  defined by the same  formulae with
$a^n,b^n,c^n,X^n(x^n,\cdot)$ replaced by $a,b,c,X(x,\cdot)$. Then
$\tilde a^n$ are uniformly Lipschitz and converge to $\tilde a$
uniformly on every compact. Thus $\sup_{r\leq
t}E\|M^n(r)-M(r)\|\to 0$. The Doob  martingale inequality provide
that $E\sup_{r\leq t}\|M^n(r)-M(r)\|\to 0$ in probability. This,
via Gronwall lemma, provides (\ref{9411}).

Next, for a given  sequences $\{n(k),k\geq1\}\subset \NN$,
$\{q(k),k\geq1\}\subset (0,+\infty)$ and $\eps>0$ we consider the
random sets
$$D_{x,t}^{*,\eps}=
\bigcup_k\{(r,v,p): r\in[0,t], |p-b(X^{n(k)}(x^{n(k)},r-),v)|\leq
\eps\theta (v)\},$$ $$ D_{x,t}^{\diamondsuit,\eps}=
\bigcup_k\{(r,v,p): r\in[0,t], |p-b(X^{n(k)}(x^{n(k)},r-),v)|\leq
\eps\theta (v),\theta(v)> q(k)\beta(v) \} $$
$$ D_{x,t}^{k,\eps}=
\{(r,v,p): r\in[0,t], |p-b(X^{n(k)}(x^{n(k)},r-),v)|\leq
\eps\theta (v),\theta(v)\leq q(k)\beta(v) \},\quad k\geq 1$$ and
put $\Omega_{x,t}^{*,\eps}=\{\forall \tau\in \Df,
(\tau,p_1(\tau),p_2(\tau))\not\in D_{x,t}^{*,\eps}\},
\Omega_{x,t}^{\diamondsuit,\eps}=\{\forall \tau\in \Df,
(\tau,p_1(\tau),p_2(\tau))\not\in D_{x,t}^{\diamondsuit,\eps}\}$,
$\Omega_{x,t}^{k,\eps}=\{\forall \tau\in \Df,
(\tau,p_1(\tau),p_2(\tau))\not\in D_{x,t}^{k,\eps}\}$. Our aim is
to construct $\{n(k)\}$ in such a way that
\be\label{9412}P(\Omega_{x,t}^{*,\eps})\to 1,\quad \eps\to 0+.\ee
Once this construction is complete, the considerations analogous
to those given in the proof of Lemma \ref{l934} can be made
uniformly over $k$ and provide (\ref{9410}). We have
$$
\Omega_{x,t}^{*,\eps}=\Omega_{x,t}^{\diamondsuit,\eps}\cap\bigcap_k
\Omega_{x,t}^{k,\eps}.
$$
Since  $\{\theta(v)\leq q\beta(v)\}\downarrow
\{\theta=0\},q\downarrow 0$ and $\theta\in L_1(\VV,\pi)$, one can
choose a monotonically decreasing sequence $\{q(k)\}$ in such a
way that $\int_{\{\theta(v)\leq
q(k)\beta(v)\}}\theta(v)\pi(dv)\leq 2^{-k},k\geq 1.$ Analogously
to (\ref{945}), we have
$$
1-P(\Omega_{x,t}^{k,\eps})\leq 1-\exp\left[-2t\eps
\int_{\{\theta(v)\leq q(k)\beta(v)\}}\theta(v)\, dv\right]\leq
2t\eps\cdot 2^{-k},\quad k\geq 1.
$$
Therefore, $P\left(\bigcap_k \Omega_{x,t}^{k,\eps}\right)\geq
1-\sum_k (2t\eps\cdot 2^{-k})=1-2t\eps$. Next, let $\eps>0$ be
fixed. By the condition $H_3^*$,
$$
\Omega_{x,t}^{2\eps}\backslash
\Omega_{x,t}^{\diamondsuit,\eps}\subset \bigcup_k\Big\{\sup_{r\leq
t}\|X^{n(k)}(x^{n(k)},r)-X(x,r)\|\geq \eps q(k)\Big\}.
$$
It follows from (\ref{9411}) that the sequence $\{n(k)\}$ can be
constructed in such a way that $P(\sup_{r\leq
t}\|X^{n(k)}(x^{n(k)},r)-X(x,r)\|\geq \eps q(k))\leq t\eps 2^{-k},
k\geq 1$. Then $$P(\Omega_{x,t}^{\diamondsuit,\eps})\geq
P(\Omega_{x,t}^{2\eps})-\sum_{k=1}^\infty P(\sup_{r\leq
t}\|X^{n(k)}(x^{n(k)},r)-X(x,r)\|\geq \eps q(k))\geq
1-\exp\left[-4t\eps \int_{\VV}\theta(v)\, dv\right]-t\eps,
$$
and therefore \be\label{9413}P(\Omega_{x,t}^{*,\eps})\geq
1-\exp\left[-4t\eps \int_{\VV}\theta(v)\, dv\right]-3t\eps.\ee We
remark that the sequence $\{n(k)\}$ has been built for a given
$\eps>0$, and for the other values of $\eps$ (\ref{9413}) may
fail. Now we proceed in the following way. We take $\eps_j=2^{-j},
j\in\NN$ and construct consequently the sequences
$\{n_1(k)\},\{n_2(k)\},\dots$ in such a way that every
$\{n_{j+1}(k)\}$ is a subsequence of $\{n_j(k)\}$ and, for every
$j$, (\ref{9413}) holds true with $\{n(k)\}=\{n_j(k)\}$ and
$\eps=\eps_j$. Then we put $\{n(k)=n_k(k), k\geq 1\}$. For this
sequence, (\ref{9413}) holds true for $\eps=\eps_j,j\in \NN$ by
the construction. This implies (\ref{9412}).  The lemma is proved.

Consider a grid $\Gf$ with $\Gamma_i=\Theta_i\times I_i,
i=1,\dots,m$, where $\Theta_1,\dots,\Theta_m\subset \VV$ are a
bounded measurable sets and $I_1,\dots,I_m\subset \ax$ are a
finite segments.

\begin{prop}\label{p332} 1. Under conditions $H_2$ -- $H_4$,  every component of the
vector $X(x,t)$ is a.s. differentiable w.r.t. $\Gf$ for every
$x\in\Re^m,t\in\ax$ and \be\label{9414}(D_i^\Gf
X_j(x,t))_{j=1}^m=\int_0^t\iint_{\Gamma_i\cap(\VV\times [0,
b(X(x,s-),v)])}\Ef_{s}^t\Delta(X(x,s-),v)Jh_i(s)\,\nu(ds,du),\quad
i=1,\dots,m.
  \ee
 2. Let conditions $H_2^*$ -- $H_4^*$, $H_5$ hold true and sequences $x^n\to x, t^n\to t$ and
 $\{n'(k), k\geq 1\}\subset \NN$ be given. Then there exist subsequence
 $\{n(k), k\geq 1\}\subset\{n'(k), k\geq 1\}$ and sets
 $\Omega_j'\in\Ff, j\geq 1$ such that $P(\bigcup_j\Omega_j')=1$
 and the
 family $\{X^{n(k)}(x^{n(k)},t^{n(k)}), k\geq 1\}$ has
 uniformly dominated increments w.r.t. $\Gf$ on every $\Omega'_j$.
\end{prop}
\demo  Let $i\in\{1,\dots,m\}$ be fixed, denote
$\nu^i=\nu(\cdot\backslash (\ax\times \Gamma_i))$.  One can
replace $\nu$ by $\nu_i$  in SDE (\ref{9415}) and apply Lemma
\ref{l934} for this new equation. Then, by the statements 1,2 of
this lemma and standard theorem on measurable modification, there
exists a function $\Psi^i:(s,t,x,\omega)\mapsto
\Psi^i_{s,t}(x,\omega)$ such that

1) $\Psi^i$ is $\Bf(\ax)\otimes\Bf(\ax)\otimes \Bf(\Re^m)\otimes
\sigma(\nu^i)$ -- $\Bf(\Re^m)$ measurable;

2) for every $s\leq t,$ the function $\Psi^i_{s,t}$ is $
\Bf(\Re^m)\otimes \sigma(\nu^i|_{[s,t]\times \VV\times \ax})$ --
$\Bf(\Re^m)$ measurable;

3) the process $X^i(x,\cdot)= \Psi^i_{s,\cdot }(x)$ is the
solution to (\ref{9415}) with $\nu$ replaced by $\nu^i$.

\noindent Define the function $\Phi:\ax\times\VV\times
\ax\times\Re^m$ by
$$
\Phi:(s,v,p,x)\mapsto \Phi_{s,v,p}(x)\eqdef x+c(x,v)\1_{p\leq
b(x,v)}.
$$
Then solution of (\ref{9311}) can be represented at the form
$$
X(x,t)=\Big[\Psi^i_{0,\tau_1^i}\circ
\Phi_{\tau_1^i,p_1(\tau_1^i),p_2(\tau_1^i)}\circ\dots\circ
\Phi_{\tau_k^i,p_1(\tau_k^i),p_2(\tau_k^i)}\circ
\Psi^i_{\tau_k^i,t}\Big](x),
$$
where $\tau_j^i\eqdef\tau_j^{\Gamma_i}$ and $k$ is such that
$\tau_k^i\leq t, \tau_{k+1}^i>t$. Let us deduce statement 1 from
Lemma \ref{l934}. In order to shorten notation, we suppose $k=1$,
the general case can be treated using the technique involving
partitions $Q\in\Qf_k$, introduced in the proof of Proposition
\ref{r41}.  The variables $\tau_j^i,p_1(\tau_j^i),p_2(\tau_j^i)$
are jointly independent with $\nu^i$. In addition, $p_2(\tau_1^i)$
possesses a distribution density, thus
$p_2(\tau_1^i)\not=b(\Psi^i_{0,\tau_1^i}(x), p_1(\tau_1^i) v)$
a.s. Therefore, since ${\prt\over \prt
s}|_{s=0}T_s^i\tau_1^i=-Jh_i(\tau_1^i)$ and  $T_s^i$ and does not
change $\nu^i, p_1(\tau_1^i)$ and $p_2(\tau_1^i)$, it is enough to
prove that, for every $s,v,p$, \be\label{9416} {\prt\over \prt
s}\Big[\Psi^i_{0,s}\circ \Phi_{s,v,p}\circ
\Psi^i_{s,t}\Big](x)=-\Ef_s^t\Delta(\Psi^i_{0,s}(x),v)\1_{\{p<b(\Psi^i_{0,s}(x),v)\}}\quad
\hbox{a.s. on the set }\{p\not=b(\Psi^i_{0,s}(x),v)\}. \ee
Statement 3 of Lemma \ref{l934} and the chain rule  yield that
\be\label{9417} {\prt\over \prt s}\Big[\Psi^i_{0,s}\circ
\Phi_{s,v,p}\Big](x)=a(\Psi^i_{0,s}(x))+\nabla
c(\Psi^i_{0,s}(x),v)a(\Psi^i_{0,s}(x))\1_{\{p<b(\Psi^i_{0,s}(x),v)\}}\quad
\hbox{a.s. on the set }\{p\not=b(\Psi^i_{0,s}(x),v)\}. \ee

Denote, by $\chi$, the distribution of $\Big[\Psi^i_{0,s}\circ
\Phi_{s,v,p}\Big](x)$. It follows from the statement 3 of Lemma
\ref{l934} and Fubini theorem that, on some set $\Omega^i_s\in
\sigma(\nu^i|_{[s,t]\times \VV\times \ax})$ with
$P(\Omega^i_s)=1$,
$${\prt\over \prt y}\Psi^i_{s,t}(y)=\Ef_s^t,\quad {\prt\over \prt
s}\Psi^i_{s,t}(y)=-\Ef_s^ta(y)\hbox{ for $\chi$-almost all $y$}.$$
Since $\Big[\Psi^i_{0,s}\circ \Phi_{s,v,p}\Big](x)$ and
$\nu^i|_{[s,t]\times \VV\times \ax}$ are independent, this
together with (\ref{9417}), Fubini theorem and the chain rule
provides (\ref{9416}). This proves statement 1.

In order to shorten notation, we prove statement 2 in the case
$m=1$. In this case, $\Gf=\{[a_1,b_1),h_1,\Gamma_1\}$. In
addition, we consider the particular case $\omega\in\tilde
\Omega=\{\#(\Df^{\Gamma^1}\cap [0,t^n])=1, n\geq 1\}$. On can see
that the considerations given below can be extended to
$\omega\in\Omega$ and arbitrary $m$ straightforwardly.
 Consider the functions
$\Psi^{1,n}$ satisfying the properties 1),2) listed above (with
$i=1$), such that the processes $X^{1,n}(x,\cdot)=
\Psi^{1,n}_{s,\cdot }(x)$ are the solutions to (\ref{9415}) with
$\nu$ replaced by $\nu^1$ and $a,b,c$ replaced by $a^n,b^n,c^n$.
Put $\Phi_{s,v,p}^n(x)\eqdef x+c^n(x,v)\1_{p\leq b^n(x,v)}.$

Denote by $\tau$ the (unique) point from $\Df^{\Gamma^1}\cap
[0,t^n]$.  Then
$$
X^n(x^n,t^n)=\Big[\Psi^{1,n}_{0,\tau}\circ
\Phi_{\tau,p_1(\tau),p_2(\tau)}^n\circ\Psi^{1,n}_{\tau,t_n}\Big](x^n).
$$
By taking a subsequence, we can suppose that $\sup_{r\leq \sup
_nt^n}\|X^n(x^n,r)-X(x,r)\|\to 0$ almost surely, see (\ref{9411}).
 With probability 1,
$\Df^{\Gamma^1}\cap \{t\}=\emptyset$ and thus a.s. there exists
$n=n(\omega)\in \NN$ and $\eps=\eps(\omega)>0$ such that
$t_n>T_{lh_1}\tau, t_n<T_{lh_1}\tau_2, |l|<\eps$, where $\tau_2$
denotes the second point from $\Df^{\Gamma^1}$. Denote
$\tau_l=T_{-lh_1}\tau, T_l=T_{lh_1}^{\Gamma_1}, l\in\Re.$ Since
$T_l$ does not change $p(\tau)$ and $\Psi^{1,n}$, we have
$$
T_l X^n(x^n,t^n)=\Big[\Psi^{1,n}_{0,\tau_l}\circ
\Phi_{\tau_l,p_1(\tau),p_2(\tau)}^n\circ\Psi^{1,n}_{\tau_l,t^n}\Big](x^n),\quad
|l|<\eps.
$$
For every $s<r$ and $y\in\Re^m$, we have
$\Big|\Psi_{s,r}^{1,n}(y)\Big|\leq \eta(r)-\eta(s)$ (the process
$\eta(\cdot)$ is given in Lemma \ref{l934}). The function
$l\mapsto \eta(\tau_l)$ is a.s. continuous at the point $0$, see
the proof of Lemma \ref{l934}. In addition,
$p_2(\tau)\not=b\Big(X(x,\tau-), p_1(\tau) \Big),
p_2(\tau)\not=b^n\Big(X^n(x^n,\tau-), p_1(\tau)\Big), n\geq 1$
a.s. Thus the variable $\eps(\omega)$ can be chosen in such a way
that $\eps>0$ a.s. and
$$
p_2(\tau)<
b^n\Big(X^n(x^n,\tau_l-),p_1(\tau)\Big)\Longleftrightarrow
p_2(\tau)< b^n\Big(X^n(x^n,\tau-),p_1(\tau)\Big),\quad |l|<\eps,
$$
Since every mapping $y\mapsto y+c^n(y,p_1(\tau)), n\geq 1$ is
Lipschitz with the constant $\gamma(p_1(\tau))$, the previous
considerations  provide that, for every $l_1<l_2$ with
$|l_{1,2}|<\eps$, \be\label{9421}
\left|\Big[\Psi^{1,n}_{0,\tau_{l_1}}\circ
\Phi_{\tau_{l_1},p_1(\tau),p_2(\tau)}^n\Psi^{1,n}_{\tau_{l_1},\tau_{l_2}}\Big](x^n)-
\Big[\Psi^{1,n}_{0,\tau_{l_2}}\circ
\Phi_{\tau_{l_2},p_1(\tau),p_2(\tau)}^n\Big](x^n) \right|\leq
$$
$$
\leq L\Big[\eta(\tau_{l_2})-\eta(\tau_{l_1}))\Big] \ee with a
random variable $L=1+\gamma(p_1(\tau))<+\infty$ a.s. Therefore,
for every $j\in\NN$, there exists an a.s. positive random variable
$\vartheta_j$ such that, for every $l_1,l_2$ with
$|l_{1,2}|<\vartheta_j(\omega)$, $$
\left|\Big[\Psi^{1,n}_{0,\tau_{l_1}}\circ
\Phi_{\tau_{l_1},p_1(\tau),p_2(\tau)}^n\circ
\Psi^{1,n}_{\tau_{l_1},\tau_{l_2}}\Big](x^n)-
\Big[\Psi^{1,n}_{0,\tau_{l_1}}\circ
\Phi_{\tau_{l_1},p_1(\tau),p_2(\tau)}^n\Big](x^n) \right|< j^{-1}.
$$

The construction of the sequence $\{n(k)\}$ at the end of the
proof of Lemma \ref{l935} can be modified easily  in order to
provide an additional requirement $\{n(k)\}\subset \{n'(k)\}$. We
suppose further that $n'(k)=n(k)=k$ (this supposition is made for
notational convenience, only). Let $\zeta,\sigma$ be the random
variables given by Lemma \ref{l935}. Denote $\Omega_j'=\{\zeta\leq
j , \sigma\geq j^{-1}\}$, then $P(\bigcup_j\Omega'_j)=1$. For
every given $l\in\Re$ and $r>0$, we have that, on the set
$\{\tau_l<r\}$,
$$
\Big[\Psi^{1,n}_{0,\tau_l}\circ
\Phi_{\tau_l,p_1(\tau),p_2(\tau)}^{n}\circ\Psi^{1,n}_{\tau_l,r}\Big](x^n)=
T_{l}\Big(X^n(x^n,r)\Big).
$$
Thus we can apply Lemma \ref{l935} and write that, for
$T_{l}\omega\in \Omega_j'$ and $r>\tau_{l}$, \be\label{9418}
\left|\Big[\Psi^{1,n}_{0,\tau_{l}}\circ
\Phi_{\tau_{l},p_1(\tau),p_2(\tau)}^n\circ\Psi^{1,n}_{\tau_{l},r}\Big](x^n)-
\Psi^{1,n}_{\tau_{l},r}(y)\right|\leq
j\left|\Big[\Psi^{1,n}_{0,\tau_{l}}\circ
\Phi_{\tau_{l},p_1(\tau),p_2(\tau)}^n\Big](x^n)-y \right| \ee as
soon as $\left |\Big[\Psi^{1,n}_{0,\tau_l}\circ
\Phi_{\tau_l,p_1(\tau),p_2(\tau)}^n\Big](x^n)-y \right|\leq
j^{-1}$. Then we apply (\ref{9418}) with $r=t,l=l_2$,
$y=\Big[\Psi^{1,n}_{0,\tau_{l_1}}\circ
\Phi_{\tau_{l_1},p_1(\tau),p_2(\tau)}^n\circ\Psi^{1,n}_{\tau_{l_1},\tau_{l_2}}\Big](x^n)$
and (\ref{9421}) in order to get
$$
\Big|T_{l_1} X^n(x^n,t^n)-T_{l_2} X^n(x^n,t^n)\Big|\leq g(l_1\vee
l_2)-g(l_1\wedge l_2),\quad |l_{1,2}|\leq \varrho, \quad
T_{l_{1,2}}\omega\in \Omega_j'\cap \tilde \Omega
$$
with $\varrho=\vartheta_j\wedge \eps, g(l)=jL\eta(\tau_l)$. The
proposition is proved.

{\it Proof of Theorem \ref{t934}.} The statement 1 can be deduced
from Theorem \ref{t212} and the first statement of Proposition
\ref{p332}; the arguments here are literally the same with those
used in the proof of statement 1 of Theorem \ref{t933}.
Analogously, Theorem \ref{t218}, the second statement of
Proposition \ref{p332} and  the arguments used in the proof of
statement 2 of Theorem \ref{t933} provide that, for every sequence
$\{n'(k)\}$ there exists a subsequence $\{n(k)\}\subset \{n'(k)\}$
such that
$$
P|_{A}\circ [X^n(x^n,t^n)]^{-1}\tov P|_{A}\circ
[X(x,t)]^{-1},\quad k\to+\infty, \quad A\subset \Nf(f).
$$
This yields the statement 2 of Theorem \ref{t934}. At last, denote
$\VV_0=\{v: \gamma(v)> {1\over 2}\}, \UU_0=\{u=(v,p): v\in\VV_0,
p\leq \beta(v)\}$ . It follows from the condition $H_3$ that
$\|\nabla c(y,v)\|\leq {1\over 2}$ for $y\in\Re^m, v\not\in
\VV_0$, and thus the matrices $\Ef_r^t, r\in[s,t]$ are invertible
as soon as $\nu([s,t]\times \UU_0)=0$ (e.g. \cite{protter},
Chapter 5, \S 10). Then we apply the arguments  used in the proof
of statement 3 of Theorem \ref{t933} on the time interval $[s,t]$
and obtain that a.s.
$$
N(f)\supset
\left\{\Span\left\{\Ef_\tau^t\Delta(X(\tau-),p_1(\tau)),\tau\in\Df\cap[s,t]:
p_2(\tau)\in\Big[0,b(X(\tau-),
p_1(\tau))\Big]\right\}=\Re^m\right\}\supset
$$
$$
\supset \Big\{\nu([s,t]\times \UU_0)=0\Big\}
$$
for every $s\in[0,t]$. We have
$$P(\nu([s,t]\times
\UU_0)=0)=1-\exp\left[(t-s)\int_{\VV_0}\int_0^{\beta(v)}dp\,
\pi(dv)\right]\geq 1-\exp\left[{t-s\over
2}\int_{\VV_0}\beta(v)\gamma(v)dp\, \pi(dv)\right]\to 0, \, s\to
t-,
$$
that proves statement 3. The theorem is proved.

\subsection{A discussion on the differential properties of the solution to SDE with
jumps.} The given above proofs of Theorems \ref{t933},\ref{t934}
are based on Theorem \ref{t218}. The  question  whether the
Theorem \ref{t218} in the proofs can be replaced  by Theorem
\ref{t217} is quite natural, since the latter theorem has a
simpler formulation and does not require any additional technical
assumptions like the one given in Definition \ref{d31}. For an SDE
with additive noise, the answer is positive: under additional
condition $\sup_{s\leq t}E\|Z(s)\|^m<+\infty$, Proposition
\ref{p33} and (\ref{943}) provide $L_m$ differentiability of every
component of $X(x,t)$ via the dominated convergence theorem. One
can use Theorem \ref{t217} in order to prove statements 1,2 of
Theorem \ref{t933} under this additional supposition and then
remove this supposition via a localization procedure.

The following example shows that, for an SDE with non-additive
noise, even when the jump rate is constant, the situation is
essentially different.

\begin{ex}\label{e41}
Suppose that $\UU=\{u_1,u_2\}$ and $\Pi(\{u_{1,2}\})\in
(0,+\infty)$ and consider SDE
$$
X(x,t)=x+\int_0^tc(X(x,s-),u)\nu(ds,du)
$$
with $c(x,u_1)=1, c(x,u_2)=\mathrm{arctg}\, x.$ One can see that
this SDE can be considered as an equation of the type (\ref{9311})
and $H_2$ -- $H_4$ hold true. Consider the grid $\Gf=\{[0,t],
h,\Gamma\}$ with $\Gamma=\{u_1\}$, let us show that $X(x,t)$ is
not $L_p$-differentiable w.r.t. $\Gf$ for any $p\geq 1$. It
follows from Lemma \ref{l32} that, for any $f$ being
$L_1$-differentiable,  the function $\Re\ni l\mapsto T_{th}^\Gamma
f(\omega)\in \Re$ belongs to $W_{1,loc}^1(\Re,\lambda^1)$ and
therefore is continuous for almost all $\omega\in\Omega$.
 Consider the set $\Omega_{1,1}=\{\nu([0,t]\times \{u_{i}\})=1,
i=1,2\}$. This set has positive probability. In addition, on this
set,
$$T_{lh}^\Gamma X(x,t)=\phi(T_{-lh}\tau_1,\tau_u),\quad
\phi(s_1,s_2)\eqdef
\begin{cases}x+1+\mathrm{arctg}\,(x+1),&s_1<s_2\\
x+\mathrm{arctg}\,x+1,&s_1>s_2\end{cases},
$$
where $\tau_1$ ($\tau_2$) denotes the time moment of unique jump
with the value of the jump equal $u_1$ ($u_2$). Thus the function
$\Re\ni l\mapsto T_{lh}^\Gamma X(x,t)$ is discontinuous with
positive probability and therefore $X(x,t)$ fails to be $L_1$
differentiable w.r.t. $\Gf$.
\end{ex}

We remark that Example \ref{e41} gives a counterexample to Theorem
2.1.3 \cite{Denis}, also.  Let,  in the notation of \cite{Denis},
$M_1=\{u_1\}, M_2=\{u_2\}, M=M_1\cup M_2$, then one has $X(x,t)\in
d$ but $X(x,t)\not \in \tilde d$. The latter statement follows
from the definition of the tensor product of two Dirichlet forms
(\cite{Denis}, section 2.1.1), Proposition 1.2.2 \cite{Denis} and
the fact that, for every $s_2\in(0,1)$,  the function
$$
[0,1]\ni s\mapsto \phi(s,s_2)
$$
is discontinuous and therefore does not belong to $H(\Delta_1)$.
This "gap" is crucial and, in general, $ (\DD^n,\Ef^n)\not\subset
(\tilde \DD^n,\Ef^n)\subset ( \DD,\Ef)$ (\cite{Denis}, subsection
2.2). Thus Theorem 3.2.2 \cite{Denis} does not imply Theorem 3.3.1
\cite{Denis}, the latter being crucial for the proof of Theorem
3.3.2 \cite{Denis}.

The situation exposed in Example \ref{e41} is quite typical, and
the solution to SDE with non-additive noise, in general, neither
is  $L_1$-differentiable in the sense of Definition \ref{d21}
above, nor belong to the domain $\Ef$ of the Dirichlet form $\DD$
in the sense of \cite{Denis}. Thus the notion of a.s. derivative
appears to be an efficient tool that allows one to consider the
functionals with quite poor differential properties and extends
 the domain of possible applications significantly.

\end{document}